\documentclass[11pt]{article}
\usepackage{amsmath, amsthm, amsfonts, amssymb}
\usepackage{url}
\newcommand{\R}{\mathbb R}
\newtheorem{theorem}{Theorem}[section]

\newtheorem{remark}[theorem]{\it \rmfamily Remark}
\newtheorem{definition}[theorem]{\it \rmfamily Definition}
\newtheorem{example}[theorem]{\it \rmfamily Examples}

\newtheorem{lemma}[theorem]{Lemma}
\newtheorem{hypothesis}{\it \rmfamily Hypothesis}
\newtheorem{proposition}[theorem]{Proposition}

\newcommand{\ble}{\begin{lemma}}
\newcommand{\ele}{\end{lemma}}
\newcommand{\bre}{\begin{remark}}
\newcommand{\ere}{\end{remark}}
\newcommand{\beq}{\begin{equation}}
\newcommand{\eeq}{\end{equation}}
\newcommand{\bth}{\begin{theorem}}
\newcommand{\bpr}{\begin{proposition}}
\newcommand{\epr}{\end{proposition}}

\def\P{{\mathbb P}}
\def\R{{\mathbb R}}

\def\C{{\mathbb C}}
\def\E{{\mathbb E }}

\def\lan{{\langle}}
\def\ran{{\rangle}}
\def\hh{{\vskip 0.5 mm \noindent }}
\def\vv{{\vskip 1mm}}

\def\qed{\hfill \hbox{\hskip 6pt\vrule
width6pt height6pt depth1pt  \hskip1pt}
\smallskip}

\begin{document}

\title {   {  Stochastic flow for SDEs with jumps and irregular drift term} }

\author{  E. Priola
  \footnote{
  %Partially supported  by the
%M.I.U.R. research project Prin 2008 ``Deterministic and stochastic
%methods in the study of evolution problems''.
E-mail:  enrico.priola@unito.it}
 \\ \\ { Dipartimento di Matematica ``Giuseppe Peano''} \\  
{Universit\`a di Torino}
 \\ { via Carlo Alberto 10,   Torino,  Italy} }

\date{}

\maketitle

%\begin{center}

%{\large Enrico Priola}\footnote
% {Dipartimento di Matematica ``Giuseppe Peano'',
%Universit\`a di Torino,
%  via Carlo Alberto 10,  \ 10123, \ Torino, Italy.
% e-mail \ \ enrico.priola@unito.it. 
% Partially supported by  the Italian National Project
% MURST ``Equazioni di Kolmogorov''.
%}

%{\vskip 8mm}

%\end{center}

{\vskip 7mm }

%\noindent {\bf Key words:} Diffusion  semigroups, degenerate
%parabolic equations, Malliavin Calculus.

\vspace{2.5 mm}

%\title{ \Large \bf Stochastic flow for SDEs with jumps and irregular %drift term\thanks{ 2010
%MSC  60H10, 34F05;
%60J75, 35B65.
%\ \thanks{
%Dedicated to Prof. J. Zabczyk  on the occasion of  50 years of %scientific activities.
%Supported by the M.I.U.R. research project Prin 2010
% } }

%\author{E. Priola
%}

%\date{}

%\maketitle

%\begin{center}
%\textbf{Abstract}
%\end{center}
%This article
%is a generalization of
% continues an earlier work of the author.
% Further applications  to existence of densities for %the solutions are %investigated.
% Finally, when $\alpha \in [1,2)$ and $\beta >0$, we  % %prove uniqueness in law.
 %$ \beta > 1 - \alpha$.
%\vspace{2.5 mm} \noindent {\bf Key words:} stochastic differential
%equations, stable processes, pathwise uniqueness, H\"older
%continuity.

\noindent {\bf Abstract:} We consider non-degenerate SDEs with a $\beta$-H\"older continuous  and bounded drift term and driven by a L\'evy noise $L$ which is of $\alpha$-stable type.
 If $\alpha \in [1,2)$ and $\beta \in ( 1 - \frac{\alpha}{2},1) $  we show pathwise uniqueness and  existence of
a  stochastic flow.
%In particular, existence of strong solutions and %pathwise uniqueness hold.
We follow the approach of [Priola,
Osaka J. Math. 2012]   improving the assumptions  on the noise
$L$.  In our previous  paper $L$ was assumed to be non-degenerate, $\alpha$-stable and symmetric. Here   we can also recover relativistic and truncated stable processes  and some classes of temperated stable processes.

\section{Introduction}

We consider the 
 SDE
\beq \label{SDE}
 X_{t} = x + \int_{0}^{t}b\left( X_{r}\right)  dr \, + \, L_{t},
 \quad
\eeq
  $x \in \R^d,$ $d \ge 1$, $ t \ge 0,$ where $b: \R^d \to \R^d $ is bounded and
 H\"older continuous of index $\beta$, $\beta \in (0,1)$,
  and $L=
 (L_t) $ is a non-degenerate
 $d$-dimensional   L\'evy process of  $\alpha$-stable type. 
 Our  main result  gives conditions under which  
  strong  uniqueness holds and, moreover, there exists a  stochastic flow. The
   present paper is a continuation of our previous work
    \cite{Pr10}, where \eqref{SDE} has been investigated assuming that $L$ is non-degenerate, $\alpha$-stable and symmetric. Here we can treat
     more general noises   like relativistic and truncated stable processes and  some classes of temperated stable processes (see the end of Section 3).
 We   follow
  the approach
  in   \cite{Pr10} showing that it works  in the present more general setting.
% indicating  the changes which are needed in % % %order to %treat the %present more general situation.

  There is an increasing   interest in
  pathwise uniqueness for SDEs when
  $b$ is  singular enough so  that
  the associated  deterministic equation \eqref{SDE}
  with $L= 0$ is not well-posed (see \cite{F} and the references therein).
   An important result in this direction was
 proved by Veretennikov in \cite{Ver} (see also \cite{Zv74}
  for $d=1$).
 He was able to prove pathwise  (or strong) uniqueness for \eqref{SDE}
 when $b : \R^d \to \R^d$ is only  bounded and measurable and
  $L $ is a  $d$-dimensional Wiener process.
 This theorem has been  extended in different directions
 (see, for instance, \cite{GM}, \cite{KR05},
 \cite{D},
  \cite{FGP},
  \cite{FF}, \cite{F}, \cite{DFPR}, \cite{proske}).

 The situation changes  when $L$ is not a Wiener process but is
  a symmetric
 $\alpha$-stable process, $\alpha \in (0,2)$.
  Indeed, when $d=1$ and $\alpha <1$, Tanaka,
  Tsuchiya and Watanabe proved in Theorem 3.2 of \cite{Tanaka}
  that even a bounded
  and $\beta$-H\"older continuous $b$ is not enough to ensure
  pathwise uniqueness  if
  $\alpha + \beta <1$.
   On the other hand,
  when $d =1$ and $\alpha \ge 1$, they showed pathwise
   uniqueness
  for any bounded and continuous drift term.

Pathwise uniqueness for equation \eqref{SDE} has been proved in \cite{Pr10}
  for $d \ge 1$, when $L$ is a
   non-degenerate, $\alpha$-stable and symmetric
    L\'evy process (cf. Chapter 3 in \cite{sato}),
   requiring that $\alpha \in [1,2)$ and
 $b$ is bounded and $\beta$-H\"older continuous
  with $\displaystyle { \beta > 1- \alpha/2}$ (see \cite{Za2} for an extension of this result when    $b$  belongs to  fractional Sobolev spaces and $\alpha >1$).
  The  approach in \cite{Pr10} 
  differs from
 the one in \cite{Tanaka} and is inspired by \cite{FGP}.
%  In \cite{Pr10} also  differentiability
% of solutions with  respect to initial conditions
%and the so-called homeomorphism  property are %established.
% We also mention \cite{Za2} where pathwise uniqueness for %\eqref{SDE} is %proved when $L$ is  non-degenerate, %$\alpha$-stable %and symmetric,  %$\alpha >1$, and  $b$  %belongs to  fractional %Sobolev spaces.
   There are  two main  examples of  L\'evy processes $L$ in the class considered  in \cite{Pr10}. The first one is the     $\alpha$-stable
 process  $L$  having the generator ${\cal L}$
 which is the fractional Laplacian $ -(- \triangle)^{\alpha/2}$, i.e.,
 \beq \label{stand0}
{\cal L} f(x) = \int_{\R^d} \big(
 f (x+z) -  f (x)
   - 1_{ \{  |z| \le 1 \} } \, \,  z \cdot D f (x)   \, \big)
   \frac{\tilde c_{\alpha,d} }{|z|^{d+ \alpha}}dz,
    \eeq
$x \in \R^d$, where $f$ is an infinitely differentiable functions
 with compact support and $Df(x)$ is the gradient of $f$ at $x \in \R^d$.  The second example is $L= (L^1_t , \ldots, L^d_t)$, where $(L^1_t)$,
 $\ldots, (L^d_t)$ are independent
 one-dimensional symmetric stable
 processes of index $\alpha$ (see \cite{basschen} and the references therein).
 In this case the    generator ${\cal L}$ is given by
\beq \label{stand1}
{\cal L}f(x) =  \sum_{k=1}^d \int_{\R} [ f (x+ s e_k) - f (x)
   - 1_{ \{  |s| \le 1 \} } \, s  \, \partial_{x_k}f (x)  ]
   \; \frac{c_{\alpha}}{|s|^{1+ \alpha}} ds.
\eeq
In the present paper we
generalize the class of L\'evy processes $L$ considered in \cite{Pr10} introducing
  Hypotheses \ref{nondeg} and \ref{nondeg1} (see Section 2). The first assumption requires the validity of some
 gradient estimates for the convolution Markov semigroup $(R_t)$ associated to $L$.
  This property  expresses the fact that $L$ must be non-degenerate to get the uniqueness result. The second assumption is an integrability condition on    the small jumps part of $L$.
 Such assumptions hold not only for the processes $L$ considered in  \cite{Pr10}. Indeed, for instance, the hypotheses  are satisfied by
    truncated stable processes (see \cite{KS} and the references therein), by some classes of temperated stable processes (see \cite{R}) and by relativistic stable processes (see, for instance, \cite{Ca} and \cite{Ry}).

    The following  is our main theorem (we recall the definitions of strong solution, pathwise uniqueness and stochastic flow in  Definition \ref{civuo}).
\bth
\label{uno} Let $L$ be a pure-jump L\'evy process satisfying Hypotheses \ref{nondeg} and \ref{nondeg1} with some $\alpha \in [1,2)$ (see Section 2).
Suppose that
 $b \in C_{b}^{\beta}\left(  \mathbb{R} ^{d} ; \R^d \right) $,
  with $\displaystyle{ \beta \in \big( 1 - \frac{\alpha}{2},1 } \big)$.
Then, we have:

\hh (i) Pathwise uniqueness holds for \eqref{SDE}, for any $x \in \R^d$.

\hh (ii) For any $x \in \R^d$, there exists a (unique) strong solution $(X_t)$ to \eqref{SDE}.

\hh (iii) There exists a stochastic flow $(X_t^x) $
   of class $C^1$.
 % (see Section 5) such that
% a modification of the strong solution which we still %denote by  $(X_t^x) $ such that,  $P$-a.s., for any $x %\in \R^d,$ trajectories
 %  $for any $x \in \R^d$, the mapping$
 %  $t \to X_t^x$ are
  %is
 % c\`adl\`ag
  %from $[0, \infty$ into $\R^d$
%  and we have:
%  with cadlag-trajectories and such
%   the solution starting at $x \in\mathbb{R}^{d}$, we %have:
\end{theorem}
Note that (iii) is stronger than (ii); condition (iii) implies  that  $P$-a.s,   for any $t \ge 0$, the mapping: $x \mapsto X_t^x$ is a differentiable
 homeomorphism from $\R^d$ onto $\R^d$ (cf. Section 3 in \cite{Ku}, Section V.7 in \cite{protter} and also \cite{Za1} for the case of  log-Lipschitz coefficients).
  Since $C_{b}^{\sigma} (  \mathbb{R} ^{d} , \R^d )
  \subset C_{b}^{\beta} (  \mathbb{R} ^{d} , \R^d ) $
 when $0 < \beta \le \sigma$,
 our main result holds  for any
  $\alpha \in [1,2)$  when  $
  \beta  \in (\frac{1}{2}, 1).$
  
  The proof of Theorem \ref{uno} is given in Section 5 and
      uses
   an It\^o-Tanaka trick (cf. Section 2 of \cite{FGP} which considers the case of a Wiener process $L$). Such method requires Schauder estimates 
  for  non-local Kolmogorov equations on  $\R^d$ like
 \beq \label{resolv}
 \lambda v(x) -
 {\cal L} v(x)  - b(x) \cdot Dv(x) = g(x), \;\; x \in \R^d,
\eeq
 with $\lambda>0$. When  $g \in C^{\beta}_b (\R^d)$, $\alpha \ge 1$
 and
 $\alpha + \beta >1$  we obtain a unique solution $v \in C^{\alpha + \beta}_b (\R^d)$ in Theorem \ref{reg}. By using suitable solutions $v$ of \eqref{resolv}, the   It\^o-Tanaka trick
 allows  to construct a diffeomorphism $\psi: \R^d \to \R^d$. This mapping $\psi$ allows
 to pass from solutions of \eqref{SDE} to solutions of an auxiliary SDE with Lipschitz coefficients (see equation \eqref{itt}). For such  equation the stochastic flow property holds.

The main difficulty of the regularity result for \eqref{resolv} is the case   $\alpha =1$. To treat such case we use a localization
     procedure which is based on Theorem \ref{sta} where Schauder
     estimates are proved in the case of  $b(x) = k$, for
     any $x \in \R^d$,
  showing that the Schauder constant is independent of $k$. Recently,
   there are many regularity results available for related non-local equations (see, for instance, \cite{basschen},
  \cite{bassJFA}, \cite{CV}, \cite{MP}, \cite{DK}, \cite{S}, \cite{MP} and the references therein); however our Theorem \ref{reg} is not covered by such results.

   It is an open problem if Theorem \ref{uno} holds even in the case  $\alpha \in (0,1)$. This is mainly due to the difficulty of proving existence of $C^{\alpha + \beta}$-solutions to \eqref{resolv} when $\alpha <1$ and $\alpha + \beta >1$. However we  mention  \cite{S} which provides $C^{\alpha + \beta}$-regularity results in the case $0< \alpha \le 1$
   for ${\cal L} = -(- \triangle)^{\alpha/2}$ using the so-called extension property of the fractional Laplacian
    (see also Remark \ref{silvan}).

\smallskip
The letters $c$ and $C$ with subscripts will
denote  positive constants 
whose values are unimportant.

\section {Assumptions and notation}

 We introduce basic concepts and notation. More details can be found
 in  \cite{A},
   \cite{sato} and \cite{Zab}.

Let  $\lan u, v \ran $ (or $u\cdot v$)  be  the euclidean inner
product between $u$ and $v \in \R^d$, for any $d \ge 1$;
  moreover
$|u|$ $ = (\lan u,u\ran)^{1/2}$. If  $C \subset \R^d$
 we denote by $1_C$ the indicator function of $C$. The
  Borel $\sigma$-algebra of $\R^d$ will be indicated by
   ${\cal B}(\R^d)$.

 Let us consider a
  stochastic basis
  $(\Omega, {\cal F}, ({\cal F}_t)_{t \ge 0}, P)$ which satisfies
  the usual assumptions (see
   \cite[page 72]{A}); the symbol  $E$  will indicate the expectation with 
respect to $P$. Recall that an
  $({\cal F}_t)$-adapted
  and
  $d$-dimensional stochastic process $L=(L_t)$ $= (L_t)_{t \ge 0}$, 
  $d \ge 1$, is a  \textit{ L\'evy  process}
  if it is  continuous in probability,
  it  has  stationary increments,
     c\`adl\`ag  trajectories, $L_t - L_s$ is independent of
     ${\cal F}_s$, $0 \le s \le t$,
      and $L_0=0$.

      One can show (see Chapter 2 in \cite{sato}) that there exists
   a unique triple $(S, b_0, \nu)$, where $S$ is  a symmetric non-negative definite $d \times d$-matrix, $b_0 \in \R^d$  and $\nu$  is a {\it L\'evy measure}
 (i.e., $\nu$ a $\sigma$-finite (positive) measure
on $\R^d$ with $ \nu (\{ 0\})=0$ and
 $  \int_{\R^d} (1 \wedge
|y|^2 ) \, \nu(dy)
 <\infty$;
 $1 \wedge |\cdot |$ $ = \min (1, |\cdot|)$) such that   the  characteristic function of
      $L_t$ verifies
 $$
\nonumber E[ e^{i \langle u, L_t \rangle }] =
 e^{- t \psi(u)} \, e^{-t \lan S u,u \ran} \, e^{i t \lan b_0, u \ran},
 $$
 \beq \label{itol}
 \psi(u)= - \int_{\R^d} \big(  e^{i \langle u,z \rangle }  - 1 -
 \, { i \langle u,z
\rangle} \, 1_{ \{ |z| \le 1 \} }  (z) \big ) \nu (dz), \;\;
 u \in \R^d, \; t \ge 0.
  \eeq
  The L\'evy measure $\nu$     is also called the
 intensity jump measure
 of $(L_t)$ and \eqref{itol} is   the
  L\'evy-Khintchine  formula.
The L\'evy measure $\nu$ of a standard $\alpha$-stable process  $L$ corresponding to  \eqref{stand0} has density $ \frac{c }{|z|^{d+ \alpha}}$; on the other hand,  the L\'evy measure  $\nu$ of  $L$ having generator in \eqref{stand1} is concentrated on axes and is singular with respect to the  $d$-dimensional Lebesgue measure.

   In this paper we only  deal with a \textit{pure-jump L\'evy process $L$ without drift term,} i.e., we assume that
\beq \label{f56}
   S=0,\;\; b_0 =0.
\eeq
  Note that, possibly changing $b(x) $ with $b(x ) + b_0$,
in equation \eqref{SDE} we may always assume that $b_0 =0$.

Thanks to \eqref{f56} we have $E[ e^{i \langle  L_t ,  u\rangle }] =
 e^{- t \psi(u)}$, $t \ge 0,$ $ u \in \R^d$; the function
  $\psi (u) $ is called the {\it symbol} (or exponent) of the L\'evy process $L$.   Given a symbol $\psi$ corresponding to a L\'evy measure $\nu$ (see \eqref{itol}) there exists a unique in law L\'evy process $M= (M_t)$ such that $E[ e^{i \langle  M_t ,  u\rangle }]$ $ =
 e^{- t \psi(u)}$, $t \ge 0,$ $u \in \R^d.$

\smallskip
    The convolution
 Markov semigroup $(R_t)$ acting on $C_b(\R^d)$ (the space of all real continuous and bounded functions on $\R^d$) and associated to
$L$ (or to $\psi$) is
\beq \label{dr66}
R_t f(x) = E [f(x+ L_t)] = \int_{\R^d} f(x + y)\,  \mu_t (dy), \;\; t>0, \; f \in
C_b(\R^d),
\eeq
 $ x \in \R^d,$ where $\mu_t$ is the law of $L_t$, and
  $R_0 = I$. Note that the Fourier transform of $\mu_t$ is $\hat \mu_t(u) $ $=\int_{\R^d}   e^{i \langle u,y \rangle }  \mu_t (dy)$ $= e^{- t \psi(u)}$, $t \ge 0,$ $ u \in \R^d.$
 The
 generator $\cal L$ of the  semigroup $(R_t)$ is given by
 \beq \label{gener}
 {\cal L}g(x) =
\int_{\R^d} \big(  g (x+y) -  g (x)
   - 1_{ \{  |y| \le 1 \} } \, \lan y , D g (x) \ran \big)
   \, \nu (dy),
 \eeq
  $g \in C^{\infty}_c(\R^d),$ where $C^{\infty}_c(\R^d)$ is the space of all
  infinitely differentiable functions with compact support
 (see \cite[Section 6.7]{A} and \cite[Section 31]{sato}); $Dg(x)$ denotes 
 the gradient of $g$ at $x \in \R^d$.

\smallskip
Before introducing the main assumptions
let us recall some  function spaces used in the paper.
 We consider  $C_{b}(\mathbb{R}
 ^{d};\mathbb{R}^{k})$, for integers
  $k,\,d\geq1$, as the set of all continuous and bounded functions
 $g:\mathbb{R}^{d}\rightarrow\mathbb{R}^{k}$. It is a Banach space endowed with the
 supremum norm $\| g\|_0 = $ $\sup_{x \in \R^d}|g(x)|,$ $g \in
  C_{b}(\mathbb{R}
 ^{d};\mathbb{R}^{k}).$
 Moreover, $C_{b}^{\beta}(\mathbb{R}
 ^{d};\mathbb{R}^{k})$, $\beta \in (0,1)$,
   is the subspace of all $\beta$-H\"older continuous
   functions $g$, i.e., $g$ verifies
$$
 [ g]_{\beta}:=\sup_{x\neq x'\in\mathbb{R}^{d}}
 \frac{|g(x)-g(x')|}{|x-x'|^{\beta}}<\infty.
$$
 $C_{b}^{\beta}(\mathbb{R}
 ^{d};\mathbb{R}^{k})$ is a Banach space with the norm
 $
  \| \cdot \|_{\beta} = \| \cdot \|_0 + [\cdot ]_{\beta}.
 $
 If $k=1$, we set
 $C_{b}^{{\beta}}({\mathbb{R}}^{d};{\mathbb{R}^k})
 = C_{b}^{{\beta}}({\mathbb{R}}^{d})$.
  Let $C_{b}^{0}(\mathbb{R}
 ^{d}, {\mathbb{R}}^{k} ) = C_{b}(\mathbb{R}
 ^{d}, {\mathbb{R}}^{k})$ and $[ \, \cdot \, ]_0 = \| \cdot\|_0$.
 For each integer $n \ge 1,$ $\beta \in [0,1)$, a function $g : \R^d \to \R$ belongs to 
 $C_{b}^{n+{\beta}}({\mathbb{R}}^{d})$
 if $g \in  C_{}^{n}({\mathbb{R}}^{d}) \cap
  C_{b}^{\beta}
 ({\mathbb{R}}^{d})$
 and, for
all $k \in \{ 1,\dots ,n\} $, the Fr\'echet derivatives $D^{k}g$
 $\in C_{b}^{{\beta}}  
 ({\mathbb{R}}^{d};{(\mathbb{R}^{d})^{\otimes (k)}} )$;
  $C_{b}^{n+{\beta}}({\mathbb{R}}^{d})$
 is a Banach space endowed with the
 norm
  $
\Vert g\Vert_{n+\beta}$ $=\Vert g\Vert_{0}+\sum_{j=1}^{n}\Vert
D^{j}g\Vert _{0}+[D^{n}g]_{\beta}$, $g\in
C_{b}^{n+{\beta}}({\mathbb{R}}^{d})$. We also define
$C^{\infty}_b(\R^d) $ $= \cap_{k \ge 1} C_b^k(\R^d).$
\begin{hypothesis} \label{nondeg}
{\em The  Markov semigroup $(R_t)$ verifies: $R_t (C_b(\R^d)) \subset
C^1_b(\R^d)$, $t>0$,
%(i.e., for any $f \in C_b(\R^d)$, $t>0$, %$R_t f$ is %continuously %differentiable and bounded),
and moreover, there exists $\alpha \in (0,2)$
and $c= c_{\alpha} >0$ (independent of $f$ and $t$), such that, for any $f \in C_b(\R^d)$,
 \begin{align} \label{grad}
 \sup_{x \in \R^d}| D R_t f(x)| \le \frac{c}{t^{1/\alpha}} \; \sup_{x \in \R^d}| f(x)|,\;\;\; t \in (0,1] .
 %\;\;\; t\in
 %(0,1),
 %f \in B_b(\R^d).
\end{align}
 }
\end{hypothesis}

\begin{hypothesis} \label{nondeg1}
{\em  For any  $\sigma > \alpha$ ($\alpha $ is the same as in \eqref{grad}), it holds
 \beq \label{tr}
 \int_{\{ |x| \le 1  \}} |x|^{\sigma} \nu (dx) < \infty.
 \eeq}
\end{hypothesis}
To prove the uniqueness result in  \cite{Pr10} it is used that the previous hypotheses
hold for non-degenerate symmetric stable processes $L$.
 However, such assumptions are satisfied in more   general cases   as the next section  shows.

\bre \label{mall} {\em  For the sake of completeness, we note that
 %we state an equivalent characterization of \eqref{grad}. %Indeed
 Hypothesis \ref{nondeg} is equivalent to the following two conditions:

\hh (i) For any $t>0$ the measure $\mu_t$ in \eqref{dr66}
 has a density $p_t$ with respect to the Lebesgue measure
  which belongs to $C^{\infty} (\R^d)$; moreover
  %$p_t(x)$ tends to 0 as $|x| \to \infty$ and
  $ |Dp_t| \in L^1(\R^d)$.

\hh (ii) We have
$$
 \int_{\R^d} |Dp_t(y)| dy \le \frac{c}{t^{{1}/{\alpha}}},\;\; t \in (0,1].
$$
It is not difficult to check that (i) and (ii) implies Hypothesis \ref{nondeg} (to this purpose one first differentiates the mapping:  $x \mapsto P_t f(x)$ when $f \in C_{c}^{\infty}(\R^d)$).

On the other hand, if we have Hypothesis \ref{nondeg} then by the semigroup and the contraction property   we deduce that, for any $t>0$, $R_t (C_b (\R^d)) \subset C^{\infty}_b(\R^d)$. Moreover, \eqref{grad} implies
$\| D^k R_t f \|_0 $ $ \le \frac{c}{(t \wedge 1)^{1/\alpha}}$ $ \| f \|_0$,  $t>0$, $k \ge 1$, $f \in C_b(\R^d)$.
 It follows that, for any $f \in C_c^{\infty}(\R^d)$, $k \ge 1$, we have
 $
\big |
\int_{\R^d} D^k f (y)  \,  \mu_t (dy ) \big|
 $ $\le c_t  \, \| f\|_0,
$
where $c_t >0$ is independent of $f$. By known  properties of the Fourier transform this implies that (i) is satisfied.

 To check (ii) we remark that, for $t \in (0,1]$, the estimate
 $$
 \big |
\int_{\R^d}  D f (y)  \,   p_t(y) dy  \big|
  =
\big |
\int_{\R^d}  f (y)  \,  D p_t(y) dy  \big|
 \le   \frac{c}{t^{{1}/{\alpha}}} \| f\|_0,
$$
which holds for any $f \in C_c^{\infty}(\R^d)$, implies (ii).

We mention that  Theorem 1 of \cite{KSc}  shows  that (i) is equivalent to   the following Hartman and Wintner condition  for the symbol $\psi$:
$$
\lim_{|y| \to \infty} \, Re  \, \psi (y) \cdot \big(\log (1 + |y|)\big)^{-1} = \infty.
$$
}
\ere

\section{Classes of L\'evy processes satisfying our assumptions}

We show examples of L\'evy processes which satisfy  our hypotheses. In particular, we concentrate on   Hypothesis \ref{nondeg}.
%. These are basically considered in Example 1.5 of %\cite{SSW}
% In \cite{P} the case of stable is considered.
%We concentrate on giving an explicit condition to get %Hypothesis \ref{nondeg}.
% We will consider the following condition:
%    {\it for $\alpha \in (0,2)$ there exists $C$, $M>0$   %such that}
%\beq\label{d34}
%Re\, \psi (y) \ge C |y|^{\alpha},\;\;\; |y| >M.
%\eeq
%To prove our result
%Gradient estimates like \eqref{grad} are investigated in %\cite{SSW}.
%In \cite{P} the case of stable is considered.
%We start with a simple  result.
\bpr \label{ciao}  Suppose that the symbol $\psi$ of $L$ (see (\ref{itol})) can be decomposed as
   $ \psi(u) $ $= \psi_1(u) + \psi_2(u)$, $u \in \R^d,$
where $\psi_1$ and $\psi_2$ are both symbols of L\'evy processes and, moreover, the  convolution Markov semigroup
 associated  to $\psi_1$ (see \eqref{dr66}$)$
%(unique in law) L\'evy process  corresponding to  %$\psi_1$
satisfies gradient estimates \eqref{grad}.

Then Hypothesis \ref{nondeg} holds.
\epr
\begin{proof} Let $t \in (0,1]$. According to
\cite[Section 8]{sato}, there exist unique infinitely divisible Borel
probability measures $\gamma_t^{(1)}$ and $\gamma_t^{(2)}$ on $\R^d$ such that the Fourier
transform
$
\hat \gamma_t^{(j)} (z) $ $  = e^{- t\psi_j (z)},$ $ j =1,2.
$
Moreover, by  \cite[Proposition
2.5]{sato} we infer that
 $ \widehat {\gamma_t^{(1)} * \gamma_t^{(2)}}$ $= \hat \gamma_t^{(1)} \, \cdot \hat \gamma_t^{(2)}$ $= e^{-t \psi}$.
 By the inversion formula
 we deduce  that $ \mu_t = \gamma_t^{(1)} * \gamma_t^{(2)}$ and so \eqref{dr66} can be rewritten as
 $$
 R_t f(x) = \int_{\R^d} \gamma_t^{(1)}(dy) \int_{\R^d} f(x + y +z)\, \gamma_t^{(2)}(dz),
 $$
 $ f \in
C_b(\R^d),
 \; x \in \R^d. $ Equivalently,
  $R_t f(x)$  $= \int_{\R^d} g_t (x+y)
    \gamma_t^{(1)}(dy)$, where
     $g_t(x)=
     \int_{\R^d} f(x  +z)\, \gamma_t^{(2)}(dz)
 $. By our assumption on $\psi_{1} $ it follows that $R_t f \in C^1_b(\R^d)$ and moreover
 $$
 |DR_t f(x)| \le \frac{c}{t^{1/\alpha}} \| g_t\|_0
  \le \frac{c}{t^{1/\alpha}} \| f\|_0,
 $$
 where $c$ is independent of $t$, $x \in \R^d$ and $f \in C_b(\R^d)$. This proves the assertion.
\end{proof}
  The next   result  follows from  Theorem 1.3 in \cite{SSW}.
\bth \label{ssw1} If for $\alpha \in (0,2)$ there exists $c_1$, $c_2$  and  $M>0$   such that
\beq \label{df44}
c_1 |y|^{\alpha} \le Re\, \psi (y) \le c_2 |y|^{\alpha},\;\;\;\; |y|>M,
\eeq
 then Hypothesis \ref{nondeg} holds.
  %for some $c_{\alpha}$.
\end{theorem}
 Condition \eqref{df44}  concerns the ``small jump part'' of the L\'evy process $L$. Indeed
if $\psi^{(1)}(u) = - \int_{\{ |y|\le 1 \}} \big(  e^{i \langle u,y \rangle }  - 1 - \, { i \langle u,y
\rangle}  \big ) \nu (dy)$ and $\psi^{(2)}= \psi - \psi^{(1)}$, then $\psi^{(2)}$ is a bounded function on $\R^d$ and so \eqref{df44} holds for $\psi$ if and only if it holds for $\psi^{(1)}$.

Using Theorem \ref{ssw1} and Proposition \ref{ciao}
(cf. Remark 1.2 of \cite{SSW}) one can obtain the following generalization of the previous result:
\bpr \label{pri1}
Assume that the L\'evy measure $\nu$ in \eqref{itol} verifies:
$
 \nu (A) \ge \nu_1(A), \;\; A \in {\cal B}(\R^d),
$
where $\nu_1$ is a L\'evy measure on $\R^d$ such that its corresponding symbol
$$
\psi_1 (u) = - \int_{\R^d} \big(  e^{i \langle u,y \rangle }  - 1 - \, { i \langle u,y
\rangle} \, 1_{ \{ |y| \le 1 \} } \, (y) \big ) \nu_1 (dy),
$$
verifies \eqref{df44}. Then Hypothesis \ref{nondeg} holds.
\epr
\begin{proof} Because the  measure $\nu_2 =\nu - \nu_1$ is still a L\'evy measure one can consider its corresponding symbol $\psi_2 = \psi- \psi_1$. Applying Proposition \ref{ciao} and Theorem \ref{ssw1} we get the assertion.
\end{proof}

\begin{example} {\em As in Example 1.5 of \cite{SSW} let $\mu$ be a finite non-negative measure on ${\cal B}(\R^d)$ with  support on the unit sphere $S$  and suppose that $\mu$ is non-degenerate (i.e., its support is not contained in a proper linear subspace of $\R^d$). Let $r>0$ and  define, for $A \in {\cal B}(\R^d)$,
\begin{gather} \label{t66}
\tilde \nu (A) = \int_0^r  \frac{ds}{s^{1+ \alpha}} \int_{S} 1_{A} (s \xi)  \mu (d \xi).
\end{gather}
It is not difficult to check that $\tilde \nu$ is a L\'evy measure on $\R^d$. Moreover $\tilde \nu$ verifies Hypothesis \ref{nondeg1}. Indeed if $\sigma > \alpha $ we have
$
\int_{\{ |x| \le 1  \}} |x|^{\sigma} \tilde \nu (dx) $ $\le
\frac{1}{\sigma - \alpha} \int_{S} |\xi|^{\sigma} \mu (d \xi)$ $
 <\infty.
$
In addition   the corresponding symbol $\tilde \psi$ verifies the assumptions of Theorem \ref{ssw1} and so Hypothesis \ref{nondeg} holds for the convolution Markov semigroup associate to $\tilde \psi$. Thus Hypotheses \ref{nondeg} and \ref{nondeg1} hold in this case.

\smallskip More generally, applying Proposition \ref{pri1}, we obtain that Hypotheses \ref{nondeg} holds if the L\'evy measure $\nu$ of the process $L$ verifies
$$
    \nu (A) \ge \tilde \nu (A),\;\;\; A \in {\cal B}(\R^d).
$$
 According to the previous discussion {\it the next examples of L\'evy processes verify Hypotheses \ref{nondeg} and \ref{nondeg1}.}

\hh - {\it $L $ is a non-degenerate symmetric $\alpha$-stable process} (see, for instance, \cite{sato} and \cite{Pr10}).

In this case $ \nu (A) = \int_0^{\infty}  \frac{ds}{s^{1+ \alpha}} \int_{S} 1_{A} (s \xi)  \mu (d \xi)$, $A \in {\cal B}(\R^d)$, $\alpha \in (0,2)$, where $\mu$ is a  finite, symmetric  measure with the support on the unit sphere $S$ which  is non-degenerate.

\hh - {\it $L $ is a truncated stable process} (see \cite{KS} and the references therein).
%(see %\cite{KS})non-degenerate symmetric $\alpha$-stable process} %(cf. \cite{Pr10}).

In this case
$$ \nu (A) = c \int_{\{ |x| \le 1\}}  \frac{1_{A}(x)}{|x|^{d+ \alpha}} \, dx,\;\;\; A \in {\cal B}(\R^d), \; \alpha \in (0,2).
$$
Note that this L\'evy measure is as $\tilde \nu$ in \eqref{t66} with $r=1$ and $\mu$ which is the normalized surface measure on $S$.

\hh - {\it $L $ is a temperated stable  process of special form} (cf. \cite{R}).
%(see %\cite{KS})non-degenerate symmetric $\alpha$-stable process} %(cf. %\cite{Pr10}).

We  consider
$$ \nu (A) =  \int_0^{\infty}  \frac{e^{-s} ds}{s^{1+ \alpha}} \int_{S} 1_{A} (s \xi)  \mu (d \xi), \;\; A \in {\cal B}(\R^d),
$$
where $\mu$ is as in \eqref{t66}, $\alpha \in (0,2)$. Note that $\nu (A) \ge \tilde \nu (A)$, $A \in {\cal B}(\R^d)$, where $\tilde \nu$ is given in \eqref{t66} with $r=1$.

\hh  - {\it $L $ is a relativistic stable process} (cf. \cite{Ca}, \cite{Ry} and the references therein).

Here $\psi (u) = \big( |u|^2 + m^{\frac{2}{\alpha}} \big)^{\frac{\alpha}{2}} - m$,  for some $m>0$, $\alpha \in (0,2)$, $u \in \R^d$. By Theorem \ref{ssw1} it is easy to see that Hypothesis \ref{nondeg} holds.  Moreover also Hypothesis \ref{nondeg1} is satisfied (see Lemma 2 in \cite{Ry}).

}
\end{example}

\section {Analytic  results for  the associated Kolmogorov equation   }

 Here   we establish  existence of regular solutions
  to equation \eqref{resolv}. This  will be
   done  through  Schauder estimates. Such regular solutions
    will be used
   to prove uniqueness for the SDE \eqref{SDE} in Section 5.

\vv {\sl It is important to remark that  Hypothesis \ref{nondeg1}
 implies that  ${\cal L}
 g (x)$
 in \eqref{gener} is  well defined  for any
 $g \in C^{1+ \gamma}_b (\R^d)$ if
 $1+ \gamma > \alpha$, $\gamma \ge 0$.  }

 Indeed ${\cal L}g(x)$
 can be decomposed into the sum of two integrals, over
 $\{ |y| >1\}$ and over $\{ |y| \le 1\}$ respectively.
 The  first integral is finite since  $g$ is bounded.
 To treat  the
 second one, we can use the estimate
 \begin{align} \label{rit}
 & | g(y + x) - g(x)
   -  \, y \cdot D g (x) |
\\ \nonumber
& \le \int_0^1 |D g (x + ry) -  D g (x)
 |\, |y| dr \le  [ Dg]_{\gamma} \, |y|^{1+  \gamma},
  \;\; |y| \le 1.
\end{align}
  In addition
  $ {\cal L} g
 \in C_b (\R^d)$ if $g \in  C^{1+ \gamma}_b (\R^d)$ and
 $1+ \gamma > \alpha$.

 We need the following   maximum principle (the proof is the
 same as  in Proposition 3.2 of \cite{Pr10}). We have to assume only  Hypothesis \ref{nondeg1} (see the discussion above).

\bpr \label{max}   Assume Hypothesis
 \ref{nondeg1} and consider $b \in C_b (\R^d, \R^d)$. If
  $v \in C^{1+ \gamma}_b (\R^d)$,  with $1+ \gamma > \alpha ,$
  $\gamma \ge 0$, is a solution
  to $\lambda v - {\cal L} v$ $- b \cdot Dv = g$, with $\lambda >0$
   and $g \in C_b(\R^d)$, then
 \beq \label{max1}
 \| v \|_0 \le \frac{1}{ \lambda} \| g
 \|_0,\;\;\; \lambda>0.
\eeq
\epr

 Next we prove  H\"older regularity  for
  \eqref{resolv} when  $b$ is  constant, i.e., $b(x)=k$, $x \in \R^d$. The general case of $b$
  H\"older continuous
   will be treated in Theorem
  \ref{reg}.
 {\sl  We stress  that the constant $c$ in \eqref{schaud}
   is independent of $b=k$.}

We impose the  natural  condition $\alpha + \beta >1$ which
 is needed to  get a regular $C^1$-solution $v$.
  On the other hand, the
  assumption $\alpha + \beta < 2 $ is not necessary in
  the next result. This condition
  simplifies the proof
 and it is not restrictive in the study
  of
    uniqueness for \eqref{SDE}.  Indeed since
  $C_{b}^{\sigma} (  \mathbb{R} ^{d} , \R^d )
  \subset C_{b}^{\beta} (  \mathbb{R} ^{d} , \R^d ) $
  when $0 < \beta \le \sigma$, it is enough to study
 uniqueness when $\beta$ satisfies  $\beta + \alpha <2 $.

\bth \label{sta} Assume Hypotheses \ref{nondeg} and \ref{nondeg1}.
 Let 
  $\beta  \in (0, 1)  $ with 
  $\alpha + \beta \in (1,2)$.
   Then,
   for any
 $\lambda> 0$, $k \in \R^d$,   $f \in C^{\beta}_b (\R^d)$,
   there exists a unique   solution
  $v= v_{\lambda} \in C^{\alpha + \beta}_b(\R^d)$ to the equation
   \beq \label{we}
 \lambda v -  {\cal L} v -
  k \cdot Dv   = f
 \eeq
 on $\R^d$.
   In addition, for any $\omega >0$ there exists $c= c(\omega)>0$ independent of
  $f$, $v$ and  $k $  such that
\beq \label{schaud}
 \, \lambda^{\frac{\alpha + \beta - 1}{\alpha }}
 \| Dv \|_{0} \, + \, [Dv]_{\alpha + \beta
-1}\le c \| f \|_{{\beta}},\;\;\; \lambda \ge \omega.
\eeq
\end{theorem}
 \begin{proof}
If $v \in C^{\alpha + \beta}_b (\R^d)$
with $\alpha + \beta>1$ then equation \eqref{we} has a meaning 
thanks to \eqref{rit}.
 Moreover,   uniqueness of solutions is a consequence of Proposition \ref{max}.

 The proof basically follows the one of Theorem 3.3
 in \cite{Pr10}. We only indicate some changes. To this purpose note that in Theorem 3.3 of \cite{Pr10} we have \eqref{schaud} for any $\lambda>0$ with $c$ independent of $\lambda$ since for the L\'evy processes  considered in \cite{Pr10} gradient estimates \eqref{grad} hold for any $t>0$ (not only for $t \in (0,1]$).

   We  first consider  the
   Markov semigroup $(P_t)$ acting on $C_b(\R^d)$ and having
${\cal L} + k \cdot D $ as generator, i.e.,
$$
P_t g(x) = \int_{\R^d} g(x \,+ y + t\, k)\,  \mu_t  ( dy ), \;\; t>0, \; \; g \in
C_b(\R^d),
$$
 $x \in \R^d,$  where $\mu_t$ is the law of $L_t$, and
  $P_0 = I$.  Then we introduce
 $v = v_{\lambda} \in C_b (\R^d)$, $\lambda >0$,
 \beq \label{d}
 v(x) = \int_0^{\infty} e^{ -\lambda \, t} P_t f( x) \, dt, \;\; x \in \R^d.
 \eeq
  We  will prove    that $v$ belongs to
    $C^{\alpha + \beta}_b (\R^d)$ and that 
\eqref{schaud} holds. Finally we will show that $v$ solves \eqref{we}.

\hh {\it I Part.} Using Hypothesis \ref{nondeg} we prove that
$v \in C^{\alpha + \beta}_b (\R^d)$
and that \eqref{schaud} holds.
  Note that $\lambda \| v\|_0 \le \| f\|_0$, $\lambda >0,$
  since each $P_t$ is a linear contraction.
 Then, we remark that
 $$
P_t g (x) = R_t \big(g(\cdot + t k ) \big)(x),  \;\; t>0,\; g \in C_b(\R^d),\; x \in \R^d,
$$
 where $(R_t)$ is defined in \eqref{dr66}. Using Hypothesis \ref{nondeg} we obtain for $t \in (0,1]$,
\beq \label{se6}
\sup_{x \in \R^d}|DP_t g(x)| \le  \frac{c \| g(\cdot + kt)\|_0} {t^{1/\alpha}}
= \frac{c \| g\|_0} {t^{1/\alpha}}.
\eeq
By using the semigroup  and the contraction property of $(P_t) $ we get easily
\beq \label{lar}
\sup_{x \in \R^d}|DP_t g(x)| \le
 \frac{c \| g\|_0} {(t \wedge 1) ^{1/\alpha}}, \;\; t>0,\;\;
 g \in C_b(\R^d).
\eeq
  Now interpolation theory ensures that   $
 \big ( C^{}_b (\R^d) \, , \,  C^{1}_b (\R^d) \big)_{\beta, \infty}
 = C^{\beta}_b (\R^d), $ $\beta \in (0,1)$,
 see, for instance, Chapter 1 in \cite{L};
   interpolating
  estimate \eqref{lar}  with the   estimate
 $\| DP_t g\|_0 \le  \| D g \|_0$, $t \ge 0,$
  $g \in C^1_b(\R^d)$, we  obtain
 \beq \label{est1}
 \| DP_t g \|_0 \le \frac{c_1}{(t\wedge 1)^{
  (1- \beta) / \alpha}}
\| g \|_{{\beta}}, \;\; t>0,\;\; g \in C^{\beta}_b(\R^d),
 \eeq
 with $c_1= c_1(c_0, \beta)$. Similarly, we  get
\beq \label{est2}
 \| D^2 P_t g \|_0 \le \frac{c_2}{(t\wedge 1) ^{
  (2- \beta) / \alpha}} \| g \|_{{\beta}}, \;\; t>0,\;\; g \in C^{\beta}_b(\R^d).
 \eeq
  By \eqref{est1} (recall that 
   $\frac{1- \beta}{\alpha} <1$)  differentiating under the
   integral sign in \eqref{d} one can easily verify 
 that $v = v_{\lambda}$ is differentiable on $\R^d$, for $\lambda >0$. 
  Moreover $Dv$ 
  is bounded on $\R^d$ and, for any $\omega >0$, there exists
  $C_{\omega}$ such that, for any $\lambda \ge \omega$ with
  $C_{\omega}>0$ independent of $v$, $k$ and $f$,
  $$
\lambda^{\frac{\alpha + \beta - 1}{\alpha }}
 \| Dv \|_{0} \le
 C_{\omega} \| f \|_{{\beta}}
  $$
   (we have used that
 $\int_0^{\infty} e^{-\lambda t} (1 \wedge t)^{-\sigma} dt $  $=
\frac{c}{\lambda^{1 -\sigma}} + \frac{e^{-\lambda}}{\lambda}$, for $\sigma < 1$ and $\lambda>0$).

\smallskip
  Finally we have to show that $Dv \in C^{\theta}_b(\R^d, \R^d)$,
  where
  $\theta = \alpha - 1 + \beta  \in (0,1)$. To this purpose  we
   proceed as in the proof
   of \cite[Proposition 4.2]{bassJFA}
    and
  \cite[Theorem 4.2]{P}.
  Using \eqref{est1},
 we find, for any $y, y' \in \R^d$, $y \not = y',$ $|y-y'| \le 1/2$,
\begin{gather*}
 |Dv (y) - Dv(y')| \le
    \int_0^{|y-y'|^{\alpha}}
  \frac{c   \| f\|_{{\beta}}}{t^{
  (1- \beta) / \alpha}} dt +
 \int_{|y-y'|^{\alpha}}^{1}
  |DP_t f(y) - DP_tf (y')| dt
 \\
 + \int_{1}^{\infty}  e^{-\lambda t}
    |DP_t f(y) - DP_t f(y')| dt.
\end{gather*}
Now using \eqref{est2} we find
$$
 |Dv (y) - Dv(y')| \le C  \| f\|_{{\beta}}
   \Big ( |y-y'|^{\theta}
    +
 \int_{|y-y'|^{\alpha}}^{1}
  \frac{|y-y'|}{t^{
  (2- \beta) / \alpha}} dt
 + \, \frac{ e^{-\lambda}} {\lambda}
   |y-y'|
 \Big)$$$$  \le c_3(\lambda) \,  \| f\|_{{\beta}}
 {|y-y'|^{\theta}},
$$
 and so $[Dv]_{ \alpha - 1 + \beta }
 \le c_3(\lambda) \,   \| f\|_{{\beta}}$, $\lambda >0$,
  where $c_3 $ is independent of  $f$, $v$,
   and $k $. Finally to get \eqref{schaud}, note that $c_3(\lambda)$ is decreasing in $\lambda$.

\vskip 1mm  \noindent  {\it II Part.}  We show that $v$ solves
  \eqref{we}, for any $\lambda >0$.

  This part can be proved  as  II Part in the proof of Theorem 3.3 in \cite{Pr10} without changes. The proof is complete.
  \end{proof}

Next we generalize   Theorem \ref{sta}
 to the case of $b \in C_b^{\beta} (\R^d, \R^d)$. 
  We can only do this
  when $\alpha \ge 1$ (cf.
   Remark \ref{al}).
 To treat the critical case $\alpha =1$ we  use a 
localization procedure.       
    This method works  since the constant
appearing  in estimate \eqref{schaud}  is independent
 of $k \in \R^d$.

\bth \label{reg} Assume Hypotheses \ref{nondeg} and \ref{nondeg1}.
  Let $\alpha \ge 1$  and
 $\beta  \in (0, 1)  $
 be such that $ \alpha + \beta \in (1,2) $.
 Then,  for any
 $\lambda>0$, $f \in C^{\beta}_b (\R^d)$,
   there exists a unique  solution
  $v= v_{\lambda} \in C^{\alpha + \beta}_b
  (\R^d)$ to 
  \beq \label{wee}
 \lambda v -  {\cal L} v  - b \cdot Dv = f
 \eeq
 on $\R^d$. Moreover, for any  $\omega >0$,
  there exists  $c = c(\omega) $, independent
  of  $f$ and  $v$,  such that
 \beq \label{sch4}
 \lambda \| v\|_0 + [Dv]_{\alpha + \beta - 1}
  \le c \| f\|_{\beta}, \;\;  \lambda \ge \omega.
 \eeq
 Finally,  we have
 $\lim_{ \lambda \to \infty} \| Dv_{\lambda} \|_{0} =0$.
\end{theorem}
 \begin{proof} The proof is the same as the proof of Theorem 3.4 in \cite{Pr10} and uses Theorem \ref{sta}.
\end{proof}

\bre \label{al}
 {\em Differently with respect to  Theorem \ref{sta},
 in Theorem \ref{reg} we are not able to prove that there exist 
 $C^{\alpha + \beta}_b$-solutions
  to \eqref{wee} when $\alpha <1$.
  This problem    is clear from
 the following simple a-priori estimate:
  $
 [ D v ]_{\alpha + \beta -1} \le
    C \| f \|_{\beta}  + C \| b \|_{\beta}  \| Dv \|_{0}
  $ $  + C \| b \|_{0}
    [ Dv ]_{\beta}.
$
Since $\alpha <1$ we only have   $Dv \in  C^{\theta}_b$
 with $\theta =  \alpha + \beta -1 < \beta $ and we cannot continue with the usual analytic methods.
 }
\ere

\section{The main  result}

 We  first  recall  basic facts and notations about
 Poisson random measures which we will use in the sequel  (see also
\cite{A}, \cite{Ku},  \cite{Zab}). We will also remind different notions  of solutions for \eqref{SDE}.

 The Poisson random measure $N$ related to 
 our L\'evy
  process $L = (L_t)$ (see  \eqref{SDE}) is given by
$$
N((0,t] \times V) \, = \, \sum_{0 < s \le t} 1_{V} (\triangle L_s) =
 \sharp \{ 0< s \le t \; : \; \triangle L_s \in V\},
$$
 for any Borel set $V$ in $\R^d \setminus \{ 0 \}$, i.e.,
  $V \in {\cal B}(\R^d \setminus \{ 0 \})$,
   $t>0$. Here
 $\triangle L _s = L_s - L_{s-}$ indicates the jump amplitude of $L$
 at  $s > 0.$
 The compensated Poisson random measure $\tilde N$ is defined by
 $$
 \tilde N ( (0,t] \times V ) \, = \,
  N((0,t] \times V) - t \, \nu (V)
$$
when $0 \not \in \bar V$ (by $ \bar V$ we denote the closure of $V$);
 recall that  $\nu $ is given in \eqref{itol}.
   By
   our assumption
  \eqref{f56} the
  L\'evy-It\^o decomposition of the process $L $
  is
 \beq \label{ito}
 L_t =
 \int_0^t \int_{\{ |z| \le 1\} } z \tilde N(ds, dz) +
 \int_0^t \int_{\{ |z| > 1 \} } z  N(ds, dz), \;\; t \ge 0.
 \eeq

 Recall that the stochastic integral
   $\int_0^t \int_{\{ |z| \le 1 \} } z \tilde N(ds, dz)$, $t \ge 0,$
    is an
  $L^2$-martingale.
 The process $ \int_0^t \int_{\{ |z|> 1\} } z N(ds, dz) $
  $   =  \sum_{0 < s \le t,
  \; |\triangle L_s| > 1}
  \triangle L _s
  $ is   a compound Poisson
process.

 Let $T>0$.
 The predictable $\sigma$-field ${\cal P}$
 on $\Omega \times [0,T]$ is generated by all  left-continuous adapted processes (defined on the same stochastic basis on which $L$ is defined).
 Let $V \in {\cal B}(\R^d \setminus \{ 0\})$ and consider   a ${\cal P} \times {\cal B}(V)$-measurable mapping  $F :
 [0, T] \times \Omega \times V \to \R^d$.

 If $0 \not \in \bar {V},$ then  $\int_0^T \int_{V} F(s,x)  N(ds, dx) $ $ = \sum_{0 < s \le T
  \; }
  F(s, \triangle L _s) 1_V (\triangle L _s)$ is a random finite sum.
If $E \int_0^T ds \int_V |F(s, x)|^2 \nu (dx) < \infty$, then one can
define
  the stochastic integral
 $$
 M_t =  \int_0^t \int_{V} F(r,x)  \tilde N(dr, dx), \;\; t \in [0,T]
 $$
 (here we do not need to assume $0 \not \in \bar {V}$).
  The process
$M = (M_t)$ is an $L^2$-martingale with a c\`adl\`ag modification. By Lemma 2.4 in \cite{Ku} we have
\beq \label{cera}E |M_t|^2
 = E\int_0^t dr \int_{V} |F(r,z)|^2  \nu(dz),\;\;\; t \in [0,T].
 \eeq

\begin{definition} \label{civuo} {\em A \textit{ weak solution} to \eqref{SDE}  is
 a tuple   $\left(
\Omega, {\mathcal F},
 ({\cal F}_t)_{t \ge 0}, P, L, X\right) $, where $(
\Omega, {\mathcal F},
 ({\cal F}_t)_{t \ge 0}, P )$ is a stochastic basis
  on which it is defined a pure-jump L\'evy process $L$ (see conditions \eqref{itol} and \eqref{f56})
  and
 a c\`adl\`ag $({\mathcal F}_{t})$-adapted $\R^d$-valued
process $X = (X_t) = (X_t)_{t \ge 0}$  which solves \eqref{SDE}  $P$-a.s..

\smallskip
A weak  solution $X$ which is $(\bar {\cal F}_t^L)$-adapted (here $(\bar{\cal F}_t^L)_{t \ge 0}$ denotes the  completed
 natural filtration of  $L$, i.e., for $t \ge 0$, $\bar{\cal F}_t^L$
 is the completed
 $\sigma$-algebra generated by $L_s$, $0 \le s \le t$)
is called \textit{strong  solution} (cf. Chapter 3 in \cite{situ}).

\smallskip
We say that \textit{pathwise uniqueness} holds for \eqref{SDE} if given two weak solutions $X$ and $Y$ (starting at  $x \in \R^d$) defined on the same stochastic basis (with respect to the same  process $L$) then, $P$-a.s., $X_t = Y_t$, for any $t \ge 0$.

\smallskip
Given a  stochastic basis $(  \Omega,
,{}\mathcal{F}, ( \mathcal{F}{}_{t})_{t \ge 0}, P )  $   on which it is defined a pure-jump L\'evy process $L$, a  \emph{stochastic flow } of class $C^1$  for \eqref{SDE}  is a map  $(t,x,\omega)\mapsto X_t^x\left(  \omega\right) $, defined
for $ t \ge 0$, $x\in{\mathbb{R}}^{d}$, $\omega\in \Omega$
with values in ${\mathbb{R}}^{d}$, such that
\begin{itemize}
\item [(a)] given $x\in{\mathbb{R}}^{d}$, the
process $X^{x}=\left( X_{t}^{x}\right)_{t \ge 0} $ is
a   c\`adl\`ag $(\bar {\mathcal F}_{t}^L$)-adapted solution of (\ref{SDE});

\item[(b)] $P$-a.s., for any $t \ge 0$, the map  $x \mapsto X_t^x$ is a
 homeomorphism from $\R^d$ onto $\R^d$;

\item[(c)] $P$-a.s.,  for any $t \ge 0$, the map  $x \mapsto X_t^x$ is
 a $C^1$-function on $\R^d$.  \qed
\end{itemize} }
\end{definition}
Starting from a stochastic flow $(X_t^x)$ one can easily obtain a stochastic flow of Kunita's type $\xi_{s,t}(x)$, $0 \le s \le t,$ $x \in \R^d$  (see Section 3.4 in \cite{Ku} and Remark 4.4 in \cite{Pr10}).

\smallskip To prove  Theorem \ref{uno} we will consider
the following   equation on
$\R^d$
\beq \label{resolv1}
 \lambda u (x)-
 {\cal L} u(x)  - Du(x) \, b(x) = b(x), \;\; x \in \R^d,
\eeq
 where $b : \R^d \to \R^d$ is given in  \eqref{SDE},
  ${\cal L}$ in \eqref{gener}
 and
 $\lambda>0$; the equation is intended componentwise, i.e., $u : \R^d \to \R^d$ and
   $$
   \lambda u_j -
 {\cal L} u_j  - b \cdot Du_j = b_j, \;\;  j =1, \ldots, d,
 $$
with $u(x) = (u_1(x), \ldots, u_d(x))$, $b(x) = (b_1(x), \ldots, b_d(x))$.

\smallskip
  The next two results only require that the drift term 
$b : \R^d \to \R^d$ is bounded and continuous.
 The first lemma provides an It\^o-type formula for a solution
 to \eqref{SDE} (cf. page 15 in \cite{FGP} where a related result is proved when $L$ is a Wiener process).

 \ble \label{due1} Let $L$ be a pure-jump L\'evy process satisfying Hypothesis \ref{nondeg1}   for some  $\alpha \in (0,2)$
 and let $b \in C_b (\R^d, \R^d)$ in  \eqref{SDE}.
 Assume
 that, for some $\lambda >0$, there exists a  solution
 $u = u_{\lambda} \in C^{1+ \gamma}_b (\R^d, \R^d)$ to
 \eqref{resolv1} with $\gamma \in [0,1]$, and moreover
 \beq
  1+ \gamma > \alpha.
\eeq
 Let $X = (X_t)$  be a weak solution of
 \eqref{SDE}, such that $X_0 =x$, $P$-a.s.. Then, $P$-a.s.,
 for any $t \ge 0$,
\begin{gather} \label{itok}
 u(X_t) - u(x)
\\
\nonumber =  x + L_t - X_t   + \lambda  \int_0^t u(X_r) dr
 +  \int_0^{t} \int_{\R^d \setminus \{ 0 \} } \! \! [ u(X_{r-} + x) - u(X_{r-})]
   \tilde N(dr, dx).
\end{gather}
\ele
 \begin{proof}  To see that  the stochastic integral in
 \eqref{itok} is well defined we  note that
   $$
  \begin{array}{l}
   E  \int_0^t dr \int_{\R^d  } | u(X_{r-} + z) - u(X_{r-})|^2  \nu(dz)
   \\ \\
  \le 4 t\| u\|_0^2   \int_ { \{ |z| >1 \}}  \nu(dz)
   +  t\| Du\|_{0}^2    \int_{ \{ |z| \le 1 \} }
   | z |^{2 } \nu(dz) < \infty.
   \end{array}
   $$
 The  assertion is obtained applying   It\^o's formula to $u_i(X_t)$, $i=1, \ldots, d$, as in  the proof of Lemma 4.2 in \cite{Pr10}
  (for more details on It\^o's formula
   see  \cite[Theorem 4.4.7]{A}
   and \cite[Section 2.3]{Ku}).

We  note that It\^o's formula (as it is usually stated) would require  that  $ u_i \in C^2 (\R^d)$. However in our situation
 we have
 \beq \label{di3}
   \int_{ \{ |x| \le 1\} }  |x|^{1+  \gamma} \nu (dx) <
\infty,
\eeq
 since $1+ \gamma > \alpha$ and  Hypothesis \ref{nondeg1} holds.  Using \eqref{di3}, \eqref{rit} and an approximation argument   one  proves that, for any  $f \in C_b^{1+\gamma}(\R^d) $,    $P$-a.s.,
  $t \ge 0,$
\begin{align} \label{ci2}
\nonumber  f(X_t) - f(x)
 &= \int_0^t \int_{ \R^d \setminus \{ 0\}} [ f(X_{s-} + z) -
f(X_{s-})]
   \, \tilde N(ds, dz)
\\ & + \int_0^t {\cal L}f(X_s)  ds + \int_0^t b(X_s)\cdot Df(X_s)  ds
 \end{align}
 (cf. It\^o's  formula (4.6) in \cite{Pr10}). Thus we can apply It\^o's formula to $u_i(X_t)$.
 
 Remark that, for any $i= $ $ 1, \ldots, d$, we have
  ${\cal L}u_i + b \cdot Du_i = \lambda u_i - b_i$ (see \eqref{resolv1}). Thus we can substitute  
  in the It\^o formula  for $u_i(X_t)$ the   term
$$
\int_0^t  {\cal L} u_i(X_r)  dr  +
 \int_0^t Du_i(X_r) \cdot b(X_r)  dr
$$
 with
$ - \int_0^t b_i(X_r) dr + \lambda  \int_0^t u_i(X_r) dr = x_i+ (L_t)_i - (X_t)_i + \lambda  \int_0^t u_i(X_r) dr
$
and obtain the assertion.
\end{proof}
 Theorem \ref{uno}
 will follow from the next  result.
\bth
\label{uno1}
 Let $L$ be a pure-jump L\'evy process satisfying Hypothesis \ref{nondeg1}   for some  $\alpha \in (0,2)$
 and let $b \in C_b (\R^d, \R^d)$ in  \eqref{SDE}.
 Suppose 
 that, for some $\lambda >0$,  there exists
 $u= u_{\lambda} \in C^{1+ \gamma}_b (\R^d, \R^d)$
 which solves  
 \eqref{resolv1} with $\gamma \in ]0,1]$, and such that
   $c_{\lambda} =
    \| Du_{\lambda} \|_{0}  < 1/3 $.
 Moreover, suppose that
 \beq
 2 \gamma > \alpha.
\eeq
Then, assertions (i), (ii) and (iii) of Theorem \ref{uno} hold.
 \end{theorem}
 \begin{proof}
% The  proof  follows the one of Theorem
%  4.3 in \cite{Pr10} with some  changes.
 We stress that  $2 \gamma > \alpha$ implies
 the condition $1 + \gamma > \alpha$ in Lemma \ref{due1}.

 Since $\|Du \|_0 <1/3 $,   the classical
  Hadamard
  theorem
  (see \cite[page 330]{protter})
  implies that
 the mapping
  $$
  \psi : \R^d \to \R^d,\;\;\;
   \psi (x) = x + u(x), \;\; x \in \R^d,
  $$
   is a
  $C^1$-diffeomorphism from $\R^d$
  onto $\R^d$.
  % (see page  the proof of Theorem
%  4.3 in \cite{Pr10}.).
 Moreover,
   $D \psi^{-1}$ is bounded on $\R^d$ with
   $\|D \psi_{}^{-1} \|_0
 \le  \frac{1}{1 - c_{\lambda}} < \frac{3}{2}$ thanks to
$$ D \psi_{}^{-1}(z)  =  [I + Du_{}( \psi_{}^{-1}
(z))]^{-1} = \sum_{k \ge 0} (- Du_{}( \psi_{}^{-1} (z)))^k,
 \; z \in \mathbb{R} ^d,
 $$
 and   we  have the estimate (see page  444 of \cite{Pr10})
\beq \label{di4}
  |D \psi^{-1}(z) - D \psi^{-1}(z')|
\le c_1  \, [ Du ]_{\gamma} \,
 |z- z'|^{\gamma}, \;\; z, \, z' \in \R^d.
\eeq
 Let $r \in (0,1]$ to be fixed later (in the proof of (i) and (ii) we may consider $r=1$) and  introduce  the following {\it auxiliary SDE} (cf. (4.11) in \cite{Pr10})
 \begin{align} \label{itt}
 & Y_t =  y + \int_0^t \tilde b ( Y_s) ds
 \\
 & \nonumber
  \int_0^t \int_{ \{ |z| \le r \} } g(Y_{s-}, z)
 \tilde N(ds, dz) +
 \int_0^t \int_{ \{ |z| > r \} } g(Y_{s-}, z)
 N(ds, dz), \;\;  t \ge 0,
\end{align}
 where
 $$
 %\begin{array}{l}
 \tilde b(y) = \lambda u( \psi^{-1} (y)) -
 \int_{ \{ |z| > r \} } \! \![ u(\psi^{-1} (y) \,  + z)  - u(\psi^{-1} (y))]
   \nu(dz)  -  \int_{ \{ r < |z| \le 1 \} } \!\! \!\! \! \! z \nu (dz)
 $$$$
   = \lambda u( \psi^{-1} (y)) -
 \int_{ \{ |z| > r \} } [g(y,z) - z ]
   \nu(dz)  -  \int_{ \{ r < |z| \le 1 \} } z \nu (dz)
 $$ and
 $$
 g(y,z) =  u(\psi^{-1} (y) \,  + z) + z - u(\psi^{-1} (y)) = \,
 $$
 $$
 = u(\psi^{-1} (y) \,  + z) + z +  \psi^{-1} (y)  -
  \psi^{-1} (y)  - u(\psi^{-1} (y))
$$
$$ = [\psi (\psi^{-1} (y) \,  + z) - y],
 \;\;\; y \in \R^d,\; z \in \R^d.
  $$
 Note that \eqref{itt}  is a SDE which satisfies the {\it usual   Lipschitz conditions} (cf.   Section 3.5 of \cite{Ku} or Section 6.2 in \cite{A}). Indeed $\tilde b $ is a Lipschitz function,  $|g(y,z)| \le $
 $(1 + \| Du\|_0) |z| $, for each $y, z \in \R^d$, and, moreover (see
  page 442 in \cite{Pr10} and Lemma 4.1 in \cite{Pr10}), for any $y, y' \in \R^d$  (recall that  $2 \gamma > \alpha$ and we are assuming \eqref{tr}),
\begin{gather*}
\int_{ \{ |z| \le 1 \}} |g(y,z) - g(y',z)|^2 \nu (dz)\le c_1 \| u\|_{1 + \gamma} 
\, |y-y'|^2 \int_{\{ |z| \le 1 \} }|z|^{2 \gamma} \nu (dz) 
\\ \le c_2 |y-y'|^2.
\end{gather*}  
(i) Let $x \in \R^d$. To prove pathwise uniqueness for our equation \eqref{SDE} note that if $(X_t)$ is a weak solution to \eqref{SDE} then
using \eqref{itok} and  formula \eqref{ito}
we easily get that $\big( \psi (X_t) \big) = \big( \psi (X_t) \big)_{t \ge 0}$ is a (strong) solution to \eqref{itt} with $y = \psi(x)$.
Since pathwise uniqueness holds for \eqref{itt}
if we consider two weak solutions  $(X_t)$ and $(Z_t)$ of \eqref{SDE} (starting at $x \in \R^d$) defined on the same stochastic basis and with respect to the same L\'evy process $L$  then we obtain, $P$-a.s., $\psi (X_t) = \psi(Z_t)$, $t \ge 0$, and so, $P$-a.s., $X_t = Z_t$, $t \ge 0.$

\smallskip \hh (ii) Let us fix  a stochastic basis on which it is defined a pure-jump L\'evy process $L$ which satisfies our  hypotheses.

Let us first consider equation \eqref{itt} and fix $r=1$. Thanks to the regularity of the  coefficients
   (cf. \cite{Ku} or \cite{A}),
for any $y \in \R^d,$   on our fixed stochastic basis there exists a  unique { strong solution} $Y^y = (Y_t^y)$.

 Let $x \in \R^d$ and set $Y_t^{\psi(x)} = Y_t$, $t \ge 0$.
  If we define $(X_t) = (\psi^{-1} (Y_t))$ we get that
  $(X_t)$ is a {\sl strong solution} to \eqref{SDE} starting at $x \in \R^d$ by  It\^o's formula.

 Standard   It\^o's formula (see Theorem 4.4.7 in \cite{A}) says that if $F \in C^2 (\R^d, \R^d)$ then we have, P-a.s., for $t \ge 0,$
 \beq \label{de5}
 \begin{array}{l}
F (Y_t) = F (\psi (x)) +
\int_0^t D F (Y_s) \, \tilde b(Y_s) ds
\\ \\
+ \int_0^t \int_{ \{ |z| \le 1 \} } \big [ F\, \big (Y_{s-} + g(Y_{s-},z)   ) - F (Y_{s-})\big]
 \tilde N(ds, dz)
\\ \\
 + \int_0^t \int_{ \{ |z| > 1 \} } \big [ F\, \big (Y_{s-} + g(Y_{s-},z)   ) - F (Y_{s-})\big]
 N(ds, dz)
\\ \\
+ \int_0^t ds\int_{ \{ |z| \le 1 \} } \big [
 F\, \big (Y_{s-} + g(Y_{s-},z)   ) - F (Y_{s-})
  - D  F (Y_s) g(Y_s, z) \big]
 \nu(dz).
\end{array}
\eeq
Recall that the third line in  \eqref{de5} is a finite random sum.  Note also that    $\int_0^t D F (Y_s) \, \tilde b(Y_s) ds $ has components
 $
\int_0^t D F_i (Y_s)  \cdot \, \tilde b(Y_s) ds, $  $i =1, \ldots$ $,d.$

\smallskip
  We fix $t>0$ and  show that  the previous formula holds even with
 $F = \psi^{-1}$ arguing by approximation (cf. page 438 in  \cite{Pr10}).
   Recall that $\psi^{-1} \in C^{1} (\R^d, \R^d)$ and $D \psi^{-1}$  is bounded on $\R^d$ and satisfies \eqref{di4}.

\smallskip
By convolution with mollifiers $(\rho_n)$, i.e.,
 $F_n (x)= \int_{\R^d} \psi^{-1}(y) \rho_n (x+y) dy$, 
  $x \in \R^d$, and, possibly passing to a subsequence, we find  $(F_n) \subset C^{\infty}(\R^d, \R^d)$ such that
  $
 F_n \to \psi^{-1}  \;\; \text{in}  \;\; C^{1+ \gamma'}(K; \R^d),
 $
  for any compact set $K \subset \R^d$ and $0 < \gamma' < \gamma$. Moreover, we have  the estimate $\| D F_n \|_0 \le  \| D \psi^{-1}\|_0$, $n \ge 1$, and,
  using \eqref{di4} and \eqref{rit}, also
\beq \label{celo}
|F_n (y+  x ) - F_n(y) - DF_n(y)x |
 \le [ D \psi^{-1}]_{\gamma} \, |x|^{1+  \gamma}, \;\;
 |x| \le 1,\; y \in \R^d,\; n \ge 1.
\eeq
Let us write It\^o' formula \eqref{de5} with $F$ replaced by  $F_n$. We can easily pass to the limit as $n \to \infty$ with $\omega \in \Omega$  fixed in the first,   third and fourth line of \eqref{de5}; for instance,
 in order to pass to the limit in the fourth line of \eqref{de5} we use estimate \eqref{celo} which gives
$$ \begin{array}{l}
  \big |
 F_n\, \big (Y_{s-} + g(Y_{s-},z)   ) - F_n (Y_{s-})
  - D  F_n (Y_s) g(Y_s, z)  \big | \le [ D \psi^{-1}]_{\gamma} \, |g(Y_s,z)|^{1+  \gamma}\\ \le C |z|^{1 + \gamma},
 \end{array}
 $$
  $|z| \le 1$, $s \in [0,t]$, $n \ge 1$,  and then we apply the Lebesgue convergence theorem taking into account   \eqref{di3}.
In the second line of  \eqref{de5} (written with $F$ replaced by  $F_n$) we  can pass to the limit in $L^2(\Omega)$ thanks to the isometry formula \eqref{cera}; indeed we have
$$
 \begin{array}{l}
E \big| \int_0^t \int_{ \{ |z| \le 1 \} }
 [ \psi^{-1} \big (Y_{s-} + g(Y_{s-},z)   )
    - F_n\, \big (Y_{s-} + g(Y_{s-},z)   )
\\ \\-   \, \psi^{-1} (Y_{s-}) + F_n \big (Y_{s-})]
 \tilde N(ds, dz) \big |^2
 \\ \\
  = E \int_0^t ds \int_{ \{ |z| \le 1 \} } \big | \psi^{-1} \big (Y_{s-} + g(Y_{s-},z)   )
    - F_n\, \big (Y_{s-} + g(Y_{s-},z)   )
\\ \\- \,  \psi^{-1} (Y_{s-}) + F_n \big (Y_{s-}) \big|^2
 \nu (dz) \to 0,
\end{array}
$$
as $n \to \infty$, thanks to the estimate
$$
\big| \psi^{-1} \big (Y_{s-} + g(Y_{s-},z)   )
    - F_n\, \big (Y_{s-} + g(Y_{s-},z)   ) - \,  \psi^{-1} (Y_{s-}) + F_n \big (Y_{s-}) \big|^2 \le C |z|^2,
$$
where $C$ is independent of $n$, $s$ and $\omega \in \Omega$ (recall that $\| DF_n \|_0 \le \| D \psi^{-1}\|_0,$ $n \ge 1$).  By It\^o's formula with $F = \psi^{-1}$ we get
\begin{gather*}
\psi^{-1}\, (Y_t) = x + \int_0^t D \psi^{-1}\, (Y_s) \tilde b(Y_s) ds
\\
+ \int_0^t \int_{ \{ |z| \le 1 \} } \big [ \psi^{-1}\, \big (Y_{s-} + [\psi (\psi^{-1} (Y_{s-}) \,  + z) - Y_{s-}]  ) - \psi^{-1} (Y_{s-})\big]
 \tilde N(ds, dz)
\\
+ \int_0^t \int_{ \{ |z| > 1 \} } z
  N(ds, dz)
\\
+ \int_0^t ds\int_{ \{ |z| \le 1 \} } \big [ z - D \psi^{-1}\, (Y_s) [\psi (\psi^{-1} (Y_{s-}) \,  + z) - Y_{s-}] \big]
 \nu(dz).
\end{gather*}
 It follows that
\begin{gather*}
 \psi^{-1}\, (Y_t) = x + L_t +  \lambda \int_0^t D \psi^{-1}\, (Y_s)   u( \psi^{-1} (Y_s))ds
\\
- \int_0^t D \psi^{-1}\, (Y_s)ds
 \int_{ \{ |z| > 1 \} } \big([\psi (\psi^{-1} (Y_{s-}) \,  + z) - Y_{s-}]  -   z \big)
   \nu(dz)
\\
+ \int_0^t ds\int_{ \{ |z| \le 1 \} } \big [ z - D \psi^{-1}\, (Y_s) [\psi (\psi^{-1} (Y_{s-}) \,  + z) - Y_{s-}] \big]
 \nu(dz).
\end{gather*}
Thus, using that $\lambda u = {\cal L} u + Du \, b + b$, we find
\begin{gather*}
\psi^{-1}\, (Y_t) = x + L_t
\\
+ \! \int_0^t D \psi^{-1}(Y_s)ds \! \!
 \int_{ \{ |z| \le  1 \} }  \! \! \! \! \big([\psi (\psi^{-1} (Y_{s-})   + z) - Y_{s-}]  -   z  - Du(\psi^{-1}\,(Y_s))z\big)
   \nu(dz)
\\
+   \int_0^t D \psi^{-1}\, (Y_s) Du (\psi^{-1}\,(Y_s)) b (\psi^{-1}\,(Y_s))ds
\\ + \int_0^t D \psi^{-1}\, (Y_s)  b (\psi^{-1}\,(Y_s))ds
\\
+ \int_0^t ds\int_{ \{ |z| \le 1 \} } \big [ z - D \psi^{-1}\, (Y_s) [\psi (\psi^{-1} (Y_{s-}) \,  + z) - Y_{s-}] \big]
 \nu(dz).
\end{gather*}
Since  $Du(y) = D \psi\, (y) - I $
 and $D \psi^{-1}(y) D \psi(\psi^{-1}(y)) = I$,  $y \in \R^d$,
we get
$$
\psi^{-1}\, (Y_t) = x + L_t + \int_0^t b (\psi^{-1}\,(Y_s))ds.
$$
 This shows that $(X_t)= (\psi^{-1}(Y_t)) $ is the unique strong solution to \eqref{SDE}.

\ \hh (iii) In order to prove the existence of a stochastic flow for \eqref{SDE}  it is enough to
  establish such property for equation \eqref{itt} and then use that
\beq \label{vivi}
  X_t^x =
  \psi^{-1}(Y_t^{\psi (x)}).
  \eeq
 To this purpose  one can
 use Theorems 3.10 and 3.4 in \cite{Ku}.
 Checking the validity of their assumptions for equation \eqref{itt}
 is quite involved and it is
 done  in the proof of Theorem 4.3 in \cite{Pr10}
 (one has also to fix $r$ in \eqref{itt} small enough; see also page 443 in \cite{Pr10} for the proof of the differentiability property).
This completes the proof.
\end{proof}

\bre \label{der} { \em We point out that the previous proof also
 provides a formula for the derivative $D X^x_t$ with respect to $x$
  in terms of $H_t^y = D Y_t^y$ (see \eqref{vivi}).
 Indeed
  $ y \mapsto Y_t^y$ is $C^1$, $P$-a.s., and
 Theorem 3.4 in \cite{Ku}  provides the following formula
  (cf. the formula at the end of page 445 in \cite{Pr10})
  $$
  \begin{array}{l}
 H_t^y  =  I + \lambda  \int_0^t Du( \psi^{-1} (Y_s^y))
  \, D\psi^{-1}(Y_s^y) \, H_s^y \, ds
 \\  + \int_0^t \int_{  \R^d \setminus \{ 0\}  }
  \big ( D_y h(Y_{s-}^y, z) \, H_{s-}^y  \big) \,
 \, \tilde N(ds, dz), \;\; t \ge 0, \; y \in \R^d.
\end{array}
$$
 According to  \cite{Pr10}
 the stochastic integral is meaningful and we have
 the estimate
  $ \sup_{0 \le s \le t} E[ |H_s|^p ] $ $< \infty$,
 for any $t>0$, $p \ge 2$.
 }
 \ere

 \noindent \textit{Proof of Theorem \ref{uno}.}
   We may assume that  $1 - \alpha/2 < \beta <  2 - \alpha$.
  We will deduce the assertion from Theorem \ref{uno1}.

  Since
  $\alpha \ge 1$, we can apply
   Theorem \ref{reg} and find
   a solution $u_{\lambda}
   \in C^{1+ \gamma}_b (\R^d, \R^d)$
   to the resolvent equation \eqref{resolv1} with
    $\gamma = \alpha - 1 + \beta \in (0,1)$.
    By
  the last assertion of Theorem \ref{reg}, we may
    choose $\lambda$ sufficiently large in order
    that $\| Du\|_0 = \| Du_{\lambda}\|_0 < 1/3$.
   The crucial assumption about $\gamma$ and $\alpha$ in
 Theorem \ref{uno1} is satisfied. Indeed
 $
 2 \gamma = 2\alpha - 2 + 2 \beta > \alpha
 $ since $\beta > 1 - \alpha/2$.
 By  Theorem \ref{uno1} we obtain the result.
 \qed

\def\ciaoo{
\bre\label{flow1}
 {\em Thanks  to Theorem
\ref{uno} we may define a stochastic flow associated to \eqref{SDE}.
  To this purpose, note that by  (ii) we have
  $X_t^x =\xi_t (x)$, $t \ge 0,$ $x \in \R^d$,
  $P$-a.s.., where $\xi_t$
  is a homeomorphism from $\R^d$ onto $\R^d$.
   Let $\xi_t^{-1}$ be the inverse map.
  As in \cite[Section 3.4]{Ku}, we set
  $
  \xi_{s,t}(x) = \xi_t \circ \xi_s^{-1}(x),\;\; 0 \le
  s \le t, \; x \in \R^d.
  $

  The family $(\xi_{s,t})$ is a stochastic flow since verifies the
 following properties ($P$-a.s):
 (i) for any $x \in \R^d$, $(\xi_{s,t}(x))$
 is a c\`adl\`ag  process with respect to $t$ and a
  c\`adl\`ag  process with respect $s$; (ii) $\xi_{s,t}:
  \R^d \to \R^d$ is an  onto homeomorphism, $s \le t$; (iii)
  $\xi_{s,t}(x)$ is  the unique solution to \eqref{SDE}
  starting from $x$ at time $s$; (iv)
   we have $\xi_{s,t}(x)
=\xi_{u,t}(\xi_{s,u}(x))$, for all $0\leq s\leq u\leq t $,
$x\in{\mathbb{R}}^{d}$, and $\xi_{s,s}(x)=x$.}
\ere
}

\bre\label{silvan} {\em  Using the  $C^{\alpha +\beta}$-regularity results
 proved in \cite{S} one can show that  the statement of Theorem \ref{reg} holds in the relevant  case of
${\mathcal L}= -(- \triangle)^{\alpha/2}$  even when $0< \alpha <1$.   Therefore,  by the   same proof given in
Section 5, one can prove that all the assertions of Theorem
\ref{uno} hold   even when $0< \alpha <1$ and $\beta
> 1 - \frac{\alpha}{2}$ if the  L\'evy process $L= (L_t)$ is a
standard symmetric and rotationally invariant $\alpha$-stable
 process.
%  With the same process $L$ also the assertion of %Corollary \ref{forte} hold with $\alpha \in (0,1)$ and % $\alpha + \beta >1.$

Nevertheless,   it remains an open problem  if
the statements of Theorem \ref{uno}
 %and Corollary \ref{forte} are
 is true when
 $0< \alpha <1$  for some other
 %in the case of  general
  non-degenerate $\alpha$-stable processes $L$
  (note that the proof of \cite{S} use the so-called extension property of $-(- \triangle)^{\alpha/2}$ which is a typical property of such non-local operator).
 For instance,
pathwise uniqueness is not clear when $\beta \in (1 - \frac{\alpha}{2}, 1)$
and $0< \alpha <1$ if $L$ has the generator as in  \eqref{stand1} and $d>1$, i.e.,
$$\mathcal L =   - \sum_{k=1}^d
   (- \partial_{x_k x_k}^2 )^{\alpha/2}.
   $$
Note that for such operator we do not know if there exist  regular solutions  $v$ to  the non-local Kolmogorov  equation \eqref{resolv}
 under the natural condition that  $b \in C^{\beta}_b(\R^d, \R^d)$, $g \in C^{\beta}_b(\R^d)$ and $\beta + \alpha >1$.
   }
\ere

\def\ciao{{\section {Introduction and Preliminaries}

Let $L$ be the Laplace operator on $\R^n$:
$$
Lu(x) = \sum_{j=1}^n   D_{jj} u(x),\;\; x \in \R^n.
$$
The classical Liouville theorem states that if for any bounded
${C^2}$-function $u$, $L u(x) =0$, $x \in \R^n$, then $u$ is
constant on $\R^n$. This result can be equivalently formulated in
terms of the heat semigroup $S_t$:
$$
S_t u(x) =  \frac{1}{ \sqrt {(2 \pi t)^n }  } \int_{\R^n} u(y)
e^{\frac{|x-y|^2 } {2t} } dy,\; t >0,  \;\; S_0 u(x) = u(x), \;x
\in \R^n.
$$
Namely  any bounded Borel function $u$, such that $S_t u(x) =
u(x)$, $ t \ge 0,$ $  x \in \R^n$, is  constant.

In more general terms, the Liouville problem can be stated as
follows. Let $E$ be a Polish space and $P_t$ be a Markovian
semigroup acting on the space ${\cal B}_b(E)$ of all real Borel
and bounded functions defined on $E$.  A map $u \in {\cal B}_b(E)
$ is said to be bounded harmonic function, shortly BHF, for $P_t$,
if
 \beq \label{he}
 P_t u(x) =u(x),\;\;\; t \ge 0, \; x \in E,
 \eeq
Under what conditions BHF for $P_t$ are constants?

Note that if $u$ is BHF for $P_t$, then  $Lu(x) =0$, $x \in E$,
where
 \beq \label{poi}
  Lu(x) = \lim_{t \to 0^+} \frac{P_t u(x) - u(x) } {t},
\;\; x \in E.
 \eeq
The operator $L$ is defined for all those functions $u \in {\cal
B}_b(E)$ for which the limit in (\ref{poi}) exists, see \cite{Dy}.
Let us recall that Markovian semigroups $P_t$ are of the form
 \beq \label{mark}
 P_t f(x) = \int_E f(y) p(t,x,dy),\;\; t \ge 0,\, x \in E,\; f \in
 {\cal B}_b (E),
 \eeq
 where $p(t,x, \Gamma)$, $t \ge 0,\, x \in E,$ $\Gamma \in {\cal B}(E)$ (${\cal B}(E)$
 being the Borel $\sigma-$field of $E$) is the  so called
 transition function. It is characterized by the properties:
\begin{eqnarray} \label{fun}
\nonumber (i ) & \; & p(t,x, \cdot) \;\; \mbox {\rm is a
probability measure on } \;\; {\cal B}(E),\;\; p (0, x, \cdot)
=\delta_x,\;\;\;  t \ge 0, \, x \in E;
\\
\nonumber (ii) & \;& p(\cdot, \cdot, \Gamma) \;\, \mbox {\rm is a
Borel function on } \;\,  [0, \infty) \times E,\;\; \Gamma \in
{\cal B}(E);
\\
\nonumber (iii) & \;& \int_E p(t,x, dy) p(s,y, \Gamma ) = p(t+s,
x, \Gamma),\;\; t,\, s \ge 0,\; \Gamma \in {\cal B}(E),
\end{eqnarray}
where $\delta_x$ is the Delta measure concentrated on $x $. The
Liouville problem has been studied for various classes of
operators $L$. The case of nondegenerate second order elliptic
operators $L$ on $\R^d$ or more generally on a differentiable
manifold $E$ has been particularly investigated, see for instance
\cite{COR}, \cite{COR1}, \cite{Cr}, \cite{Pi}, \cite{CC},
\cite{Ba}, \cite{BBG}, \cite{BL}, \cite{Wa} and references
therein.

In the present paper we investigate   Liouville problems both in
finite and in infinite dimensions, i.e. when $E = \R^n$ is
replaced by a separable Hilbert space $H$.

More precisely in Section 2, we deal with elliptic operator $L$
with diffusions operators on $\R^n$ of the form
 \begin{eqnarray} \label{lu}
\nonumber  Lu (x) & = &\frac{1}{2} \Tr (Q(x) D^2 u(x)) + \lan
F(x), D u (x)\ran
\\
\nonumber  & = & \frac{1}{2} \sum_{i, j =1}^n Q_{ij} (x) D_{ij}
u(x) +
 \sum_{i =1}^n F_i (x) D_{i}u(x), \; x \in \R^n,
 \end{eqnarray}
where $Q_{ij}$, $F_i$ are globally Lipschitz mapping on $\R^n$,
$Q(x)$ is symmetric and nonnegative definite. The Markovian
semigroup $P_t$ determinated by  $L$ is given by
$$
P_t f(x) = \E f(X^x_t),\;\; f \in {\cal B}_b(\R^n),\; t\ge 0,\, x
\in \R^n,
$$
where $X_t^x$ is the solution to the  SDE
\begin{eqnarray}                       \label{xt}
dX_t = F(X_t)dt  + B(X_t) dW_t = F(X_t)dt  + \sum_{j=1}^m B_j(X_t)
dW_t^j, \;\, X_0 = x \in \R^n,\; t \ge 0,
\end{eqnarray}
where $W_t = (W_t^j)$, $j=1, \ldots, m$, is a standard Wiener
process with values in $\R^m$, adapted with respect to a
stochastic basis $(\Omega, {\cal F}, ({\cal F}_t), \P)$ and
$B(x)B(x)^*$ $= Q(x)$, $x \in \R^n$.

The  main result of Section 2, concerns the Ornstein-Uhlenbeck
operator,   a particular case of the operator $L$ appearing in
(\ref{lu}), i.e.
$$
Lu (x) = \frac{1}{2} \Tr (Q D^2 u(x)) + \lan Ax, D u (x)\ran.
$$
%%%%Assuming that this operator is hypoelliptic
 In Theorem \ref{??} we are able to give a complete description of
 the matrices $Q$ and $A$ for which the Liouville theorem for $L$
 holds. This theorem is new even in the nondegenerate case (i.e. $Q$
 is positive definite) if $n \ge 3$ and it is one of our main
 result. Its proof is also based on our recent control theoretic
 result, see \cite{PZ}.

Section 2 contains also some nonsharp Liouville theorems  for the
more general operators $L$ as  in (\ref{lu}), see Propositions
\ref{grad}, \ref{??}.

In Sections 3, we consider an infinite dimensional
Ornstein-Uhlenbeck operator perturbed by some non-local operators
of the pseudo-differential character. Theorem \ref{??} states more
general conditions than Theorem \ref{??} implying again the
Liouville theorem. Assumptions of Theorem \ref{??} are also
motivated by the one-dimensional example ??. Let us recall that
infinite  dimensional Ornstein-Uhlenbeck processes perturbed   by
a Levy noise have recently  received a lot of attention, see for
instance \cite{Za2}, \cite{Ch}, \cite{FR}, \cite{RW} and
references therein.

In  Section  4 we show that in some generality the Liouville
theorem has an important probabilistic interpretation, see Theorem
\ref{??}. It turns out that Liouville' theorems holds true if and
only if arbitrary Borel sets are absorbing trajectories of all
Markov processes $Y_t$ corresponding to $P_t$ either with
probability 1 or with probability 0. We prove this results for
c\`adl\`ag processes $Y_t$, which take values in a Polish space
and are strong Feller and irreducible. Theorem \ref {??} is proved
in Chapter 9 of the  monograph \cite{Pi} when $L$ is a
nondegenerate elliptic operator in finite dimension, combining
analytic and probabilistic methods.

Finally Section 5 contains additional information on the subject
and gives alternative proofs to our results.

{\vskip 5mm}
  \ We denote by $| \cdot |$ and $\langle \cdot , \cdot
 \rangle $, the Euclidean norm
 and the  inner product of the  real separable Hilbert space  $H$. For any Polish
 space $E$, the space ${ \cal  C}_b ^{}(E )$ stands for the
 Banach space of all real continuous and bounded  functions \ $f:
E \to \R $, endowed with the sup norm:
 \ \ $ {
 \| f \|_0 \, =  \sup _{x \in E } | f(x) |, \;\;\; f \in
 { \cal C}_b (E ). }$

The space ${ \cal C}_b^{k }(H),$
  is the set of
  all $k$-times differentiable functions $f$, whose   derivatives
   are  continuous and bounded on $H$
 up to the order $k$. Moreover ${ \cal C}_b^{\infty }(H)$ $= \cap_{k \ge 1}$
 ${ \cal C}_b^{k }(H).$

 Let $M$ be a symmetric positive definite trace class operator,
 we denote by $N(x, M),
 $ $\; x \in H$,  the {\it Gaussian measure} on $H$
 with mean $x \in H$ and covariance operator
 $M$. Recall that if $H=\R^n$, then $N(x,M)$
 has density \ $
{1 \over \sqrt { (2 \pi)^n \, \det \, (M) } }
 \, e ^{  - \frac{1}{2}\, \langle M^{-1} (x -y), \, x - y\rangle_{}  },
$ with respect to the Lebesgue  measure $dy $.

Let $(X,  \| \cdot \|_{X })$ and  $(Y,  \| \cdot \|_{Y })$ be
 two  real Banach spaces.  ${\cal L}(X, Y)$ stands for  the
 Banach space of all bounded linear  operators from $X$ into $Y$, endowed with
 the operator   norm:  $ \| T \|_{{\cal L}(X,Y) } = \sup_ {  \| x\|_{X } \, \le 1 }
 \, \| Tx \|_{Y },\;\; T \in {\cal L}(X,Y)  $;   we set
  ${\cal L}_{}(X) = $ ${\cal L}_{}(X,X)$).

Let  $A : D(A) \subset H \to H$ be any closed operator on $H$, we
set $\sigma(A)$ to indicate the spectrum of $A$. Let $A$ be the
generator of a ${\cal C}_0$-semigroup $e^{tA}$. A closed subspace
$K \subset H$ is called  {\it invariant} for $e^{tA}$ if $e^{tA}
(K) \subset K$, $t \ge 0$. In this case we have that $A( D(A) \cap
K) \subset K$ and the restriction of $e^{tA}$ to $K$ is still a
${\cal C}_0$-semigroup on $K$ with generator $A^K: ( D(A) \cap K)
\subset K$ $\to K $, $A^K v = Av$, $v \in  D(A) \cap K $.

In the sequel we need the  following theorem which provides the
existence of useful  invariant subspaces for $e^{tA}$, see Lemma
2.5.7 in \cite{CZ} and \cite[page 245]{EN}.

Assume that $\sigma (A)= \sigma_0 \cup \sigma_1$, where
 $\sigma_0$ and $\sigma_1$ are two disjoint closed subsets of $\C$
 and in addition  $\sigma_0$ is bounded.
 Since the distance between $\sigma_0$ and $\sigma_1$ is positive,
 there exists a bounded open set $\Omega$, containing $\sigma_0$ and such that its closure
 is disjoint from $\sigma_1$. We may assume that the boundary $\gamma$ of $\Omega$ consists
 of a finite number of
  rectifiable, closed simple Jordan paths,
 oriented in the usual positive direction.
 Let us introduce the spectral Riesz projection $P_0$,
 \beq \label{pr1}
  P_0 x = \frac{1}{2 \pi i} \int_{\gamma} (w -
A)^{-1} x \, dw,\;\;\; x \in H.
 \eeq
  Define $P_0 H = E_0$
 (note that $E_0$ can be equal to $H$).
 We get  the {\it spectral decomposition}:
\beq \la{de4} H = E_0 \oplus E_1, \;\;\; E_1 = P_1 H, \;\; \mbox
{\rm where } \, \; (I- P_0)=P_1.
 \eeq
  The closed subspaces
$E_0$ and $E_1$ are invariant for $e^{tA}$ and  moreover $ E_0
\subset D(A)$.  The restrictions $A_i$ of
 $A$ to $E_i$, $i=0,1$, satisfy: $\sigma (A_i) = \sigma_i$.
   We have:
 \beq\label {ce3}
A_0 : E_0 \to E_0,\;\;\; A_1  : (D(A) \cap E_1) \subset E_1 \to
E_1.
 \eeq
 The operator $A_0$ generates a  uniformly continuous group $e^{t A_0}$ on $E_0$ and $A_1 $
 generates a ${\cal C}_0$-semigroup $e^{tA_1}$ on $E_1$. The
 restrictions of $e^{tA}$ to $E_0$ and $E_1$ coincide with $ e^{t A_0}$ and
 $e^{tA_1}$ respectively.

???????????????????

\section { Bounded harmonic functions  for finite dimensional diffusions}

We start  by  a preliminary non-sharp Liouville theorem for finite
dimensional nondegenerate diffusions. The result is not sharp as
the previous theorem but it does not follow from the results in
\cite{Pi}. It involves the so called Bismut-Elworthy-Lee formula,
see \cite{EL}, \cite{Ce}.

??? More  precisely, we are concerned with finite dimensional
processes $X_t^x$ which are solutions to
\begin{eqnarray}                       \label{n1}
dX_t = F(X_t)dt  + B(X_t) dW_t = F(X_t)dt  + \sum_{j=1}^m B_j(X_t)
dW_t^j, \;\, X_0 = x \in \R^n,\; t \ge 0,
\end{eqnarray}
where $W_t = (W_t^j)$, $j=1, \ldots, m$, is a standard Wiener
process with values in $\R^m$, adapted with respect to a
stochastic basis $(\Omega, {\cal F}, ({\cal F}_t), \P)$. Moreover
$F\in {\cal C}^2 (\R^n, \R^n)$, $B\in {\cal C}^2(\R^n, {\cal L}
(\R^m, \R^n) ) $, with first and second bounded derivatives (note
that $B_j (x) \in \R^n$, $x \in \R^n$). Under these assumptions it
is well known that there is a unique strong solution $X_t^x$ to
(\ref{n1}), see for instance \cite{IW}, \cite{Kr}, \cite{Ce}. It
is also standard to prove that
 \beq \label {est}
\E \Big ( \sup_{t \in [0,T]}  |X_t^x|^p  \Big) \le C_T(x),\;\;
T>0, \, x \in \R^n,\;\; p \ge 2.
 \eeq
 Recall the Ito formula:
\begin{eqnarray} \label{Ito}
u(X_t^x) - u(x) = \int_0^t \lan D u(X_s^x), B(X_s^x) dW_s \ran +
\int_0^t Lu(X_s) ds,\; x \in \R^n, \, t\ge 0,
\end{eqnarray}
where $u \in {\cal C}^2 (\R^n)$ and
\begin{eqnarray}                       \label{ko1}
Lu (x)= \frac{1}{2} \Tr (Q(x) D^2 u(x)) + \lan F(x), D u (x)\ran,
\; x \in \R^n;
\end{eqnarray}
here $Q(x) = B(x)B^*(x)$, $x \in \R^n$, (note that the map $Q \in
{\cal C}^2(\R^n, {\cal L} ( \R^n) ) $ and, for any $x \in \R^n$,
$Q(x) $ is nonnegative definite). The following straightforward
result clarifies the connection between BHF for $X_t$ and bounded
classical solutions of $L u=0$.
%%%% We prove it for the sake of completeness.
\begin{Proposition} \label{ciao}
Let $u \in {\cal C}^2 (\R^n) \cap {\cal C}^1_b (\R^n)$ be such that $Lu(x)=0$, $x
\in \R^n$, then u is BHF for $X_t$. Conversely if $v \in   {\cal C}^2_b (\R^n) $ is
a BHF for $X_t$, then  $Lv (x)=0$, for any  $x \in \R^n$.
\end{Proposition}
{\bf Proof} \ Assume that $u \in {\cal C}^2 (\R^n) \cap {\cal
C}^1_b (\R^n)$ be such that $Lu(x)=0$. Using (\ref{est}) we know
that the stochastic integral  $\int_0^t \lan D u(X_s^x), B(X_s^x)
dW_s \ran$ is a martingale, for any $x \in \R^n$. Taking the
expectation in the Ito formula, we  infer
$$
\E u(X_t^x) = u(x),\;\; x \in \R^n,\;\; t \ge 0.
$$
Let now $v \in   {\cal C}^2_b (\R^n) $   be  BHF for $X_t$. Taking
again the expectation in the Ito formula, and dividing by $t>0$,
we get
$$
0= \frac{1}{t} \E \int_0^t Lv(X_s^x) ds = \frac{1}{t} \int_0^t \E
Lv(X_s^x) ds, \;\;  t>0,
$$
since $v \in {\cal C}^2_b (\R^n)$. Passing to the limit as $t \to
0^+$, we obtain $Lv(x)=0$, $x \in \R^n$. \qed

Now we give a Liouville theorem for the diffusion $X_t$, assuming
 that it is non-degenerate. More precisely,  we suppose
  that $n =m$,  the matrix $B(x)$
  is invertible, for any $x \in \R^n$, and satisfies:
  \beq \label{invo}
 \sup_{x \in \R^n} \, \| B^{-1} (x) \| \; < \infty.
 \eeq
 Note that this condition is equivalent to ask that the operator
 $L$ is strongly elliptic, i.e.
 \beq \label{q1}
\inf_{x \in \R^n} \sum_{i,j =1}^n  Q_{ij} (x) h_i h_j \ge C
|h|^2,\;\; h \in \R^n,
  \eeq
for some $C>0$. Conversely recall that if we start from the
differential operator $L$,
%%% instead of $X_t$,
then, under (\ref{q1}), there always exists a map $B\in {\cal
C}^2(\R^n, {\cal L} (\R^m, \R^n) ) $, such that $B(x)^* B(x) =
Q(x)$, $x \in \R^n$, i.e. $B(x) $ is a square root of $Q(x)$.

Condition (\ref{invo}) implies  also that the process $X_t$ is
irreducible and strong Feller. More precisely,  if $P_t$ denotes
the Markov semigroup associated to $X_t$, then $P_t (I_A) (x)>0$,
for any $t>0$, $x \in \R^n$, $A $ open  set in $\R^n$, and  $P_t
({\cal B}_b (\R^n)) \subset C^2_b(\R^n)$, $t>0$.

\begin{Proposition} \label{grad}
Let $X_t$ be the solution of (\ref{n1}). Assume that (\ref{invo})
holds. Then the following condition:
 \beq \label{tesi}
\lan DF(x) k, k  \ran  \, + \,  \sum_{j=1}^m \lan DB_j (x) k ,
DB_j (x)k \ran \le 0,\; \;\; x \in R^n,\; k \in \R^n,
 \eeq
where $DB_j(x) \in {\cal L}(\R^n)$, implies that all BHF for $X_t$
are constant.
\end{Proposition}
{\bf Proof} \ Let   $P_t f (x)= \E f(X_t^x)$, $\psi \in {\cal B}_b
(\R^n)$, $t \ge 0.$ Take any BHF $\psi $ for $X_t$. Note that
$\psi \in {\cal C}_b^2(\R^n)$. Differentiating the BHF $\psi$, we
get
 \beq \label{gr3}
 \lan D\psi (x), h \ran =  \lan  D P_t \psi (x), h \ran, \;\; x,\, h
 \in \R^n.
\eeq
 Now if we prove that for any $f \in {\cal C}_b(\R^n)$, one has:
 \beq \label{gr0}
 \lim_{t \to + \infty}  \lan D P_t f (x), h \ran = 0, \;\; x,\, h
 \in \R^n,
\eeq
 then using (\ref{gr3}) we infer $\psi =$ constant. Let us verify
 (\ref{gr0}).

Since (\ref{invo}) holds, we can use  the so called
 Bismut-Elworthy-Lee formula for the first derivatives of the
semigroup $P_t$.  This reads, for any $f \in {\cal C}_b (\R^n)$,
see \cite{EL} and \cite{Ce},
\begin{eqnarray}                       \label{el1}
\langle  DP_t f(x), h\rangle = \frac{1}{t}  \E  \Big (f (X_t^x )
\int_0 ^t \langle B^{-1}(X_s^x) \,  ( [D_x X_s^x] h  ) , d W_s
\rangle \Big ), \; t > 0, \, x,\, h \in \R^n.
\end{eqnarray}
Recall that the map: $x \mapsto X_t^x$ from $x \in \R^n$ into
 $L^2 (\Omega , \R^n) $  is twice differentiable,
 thus $[D_x X(s,x)] h $ is meaningful.
Moreover we have that  $\eta (t) = [D_x X(s,x)] h$ solves the
SDE's:
$$
d \eta (t) =  D F (X_t^x)[\eta (t)] dt +   \sum_{j=1}^m D
B_j(X_t^x) [\eta (t) ] dW_t^j ,\;\; \eta (0) =h,\; t \ge 0.
$$
This is  the variation equation, obtained formally differentiating
the coefficients in (\ref{n1}). Now we apply the Ito formula to
the process $|\eta (t)|^2 $. We get
\begin{eqnarray*}
|\eta(t)|^2 - |h|^2 & =&  2 \int_0^t \lan \eta (s), DF(X_s)\ran ds
+ 2 \sum_{j=1}^n \int_0^t \lan \eta (s), DB_j (X_s) [\eta (s)]
dW^j_s \ran
\\
\; & + & \sum_{j=1}^n \int_0^t | DB_j (X_s) [\eta (s)]|^2 ds.
\end{eqnarray*}
Taking expectation and using hypothesis (\ref{tesi}), we infer
\beq \label{r5}
 \E |\eta(t)|^2 \le |h|^2,\;\; t \ge 0.
 \eeq
Using this estimate in formula (\ref{el1}), we obtain, by the
H\"older inequality,
\begin{eqnarray*}
|\langle  DP_t f(x), h \rangle|^2 \le  \frac{c \| B^{-1} \|_0}
{t^2} \| f\|_0^2 \, \int_0 ^t |h|^2 ds \le \frac{c \| B^{-1} \|_0}
{t} \| f\|_0^2 |h|^2,
\end{eqnarray*}
$t>0$, $x \in \R^n$. This gives (\ref{gr0}) and the proof is
complete. \qed

 It is worth noticing that the previous result works even in
 infinite dimensions under additional assumptions which guarantee
 the validity  of  an infinite dimensional   Bismut-Elworthy-Lee formula,
 see \cite{Za}, \cite{Ce}.

However Proposition \ref{grad} is far to be  sharp. For instance,
 if we apply this to a non-degenerate Ornstein-Uhlenbeck process,
then we get that condition (\ref{tesi}) is satisfied if $\lan Ax,x \ran$ $\le 0$, $x
\in \R^n$. This is clearly weaker than the claim contained in Theorem \ref{ot}.

Moreover it is well known that if (\ref{inv}) holds, $F $ is identically 0 and in
addition
$$
\sup_{x \in \R^n} \| B(x) \| < \infty,
$$
i.e. the operator $L$ is a purely second order uniformly elliptic operator, then the
Harnack inequlity implies that any bounded classical solution of $Lu (x)= 0$, $x \in
\R^n$, is constant, see \cite{GT}. However we can not prove this result by our
methods.

\begin{Remark}\label{grad}
{\em We note that  the method of the proof of Proposition \ref{grad},
shows more generally that under (\ref{tesi}), if for a map  $f \in {\cal B}_b(\R^n)$,
there exists $\hat t>0$ such that
$$
P_{\hat t} f(x)\, = \, \E f(X^x_{\hat t}) = f(x),\;\; x \in \R^n,
$$
then $f$ is constant. To see this, remark that the semigroup property implies that
$P_{n \hat t} f(x)$ $  = f(x)$, $n \in \N$, $x \in \R^n$.
Differentiating in $x$, we get
$$
Df(x) = DP_{n \hat t} f(x),
$$
 letting $n \to \infty$, we obtain that $f$ is constant.
Remark that this stronger form of the Liouville theorem holds
 also when the diffusion process $X_t$ is replaced by any finite dimensional
 Levy process, see \cite{CD}.
}
\end{Remark}

\subsection {A Liouville theorem for the Ornstein-Uhlenbeck operators}

Here we state  our main result in the finite dimensional case, to
be more clear, see Sections 3 and 4 for the general case. We
consider an Ornstein-Uhlenbeck process
$$
X_t ^x = e^{tA}x  + \int_0^t e^{(t-s)A} \sqrt {Q} dW_s, \;\; x \in
\R^n, \; t \ge 0,
$$
where $Q$, $A$ are matrices in $\R^n$,  $Q = Q^*$ is  non-negative
definite. Moreover $W_t$ is a standard Wiener process with values
in $\R^n$.  We assume  that $X_t$  is irreducible and strong
Feller. This is equivalent to the fact that all the matrices
$Q_t$,
 \beq \label {qtt}
 Q_t
= \int_0^t  e^{sA} Q e^{sA^*} ds, \; \; \mbox {\rm are positive
 definite}, \; \, t>0.
   \eeq
 Remark that (\ref{qtt}) is also equivalent to the Hormander condition on  hypoellipticity
of the Ornstein-Uhlenbeck differential operator $L$,
$$
Lu (x)= \frac{1}{2} \Tr (Q D^2 u(x)) + \lan Ax, D u (x)\ran, \; x
\in \R^n.
$$
It is easy to check that $u$ is a BHF for $X_t^x$ if and only if
$u$ is a classical bounded solution of $Lu(x)=0 $, $x \in \R^n$,
see Proposition \ref{ciao}.
\begin{Theorem} \label{ot}
Let us consider the Ornstein-Uhlenbeck differential operator $L$,
 assuming that (\ref{qtt}) holds.
 Then there exist
 nonconstant, bounded  classical solutions $u$ of  $Lu(x) =0$, $x \in \R^n$, if and only
 if there exists at least one eigenvalue of $A$ with positive real part.
\end{Theorem}
This result seems to be true  if $n \ge 3$. In $\R^2$, assuming
also  that $Q$ is positive  definite,
 the result follows in particular by the results in \cite{COR}. Moreover our theorem
 is also related with \cite{Me}, in which the spectrum of $L$ in spaces
 of continuous functions  is investigated.
%% (this which also describes the Martin boundary se).

\section {  Liouville theorems  in infinite dimensions
}

This section is devoted to  prove Theorem \ref{ot} for an infinite
dimensional Ornstein-Uhlenbeck process $X_t$ which solves
 the SDE's
 \beq \label{1e4} dX_t = AX_t dt + B dW_t,\;\; X_0 = x
\in H,\; t \ge 0,
 \eeq
  where $A: D(A)\subset H \to H$ generates a ${\cal
 C}_0$-semigroup $e^{tA} $ on $H$, $B \in {\cal L}(U,H)$; here $U$
is another Hilbert space. Moreover $W_t$ is an $U-$valued {\it
cylindrical} Wiener process adapted with respect to a stochastic
basis $(\Omega, {\cal F}, ({\cal F}_t), \P)$, see \cite{DZ} for
more details. This means that there exist independent one
dimensional Wiener processes $(\beta_t^k)$, $k \in \N$, such that
$W_t: U \to L^2 (\Omega )$, $W_t (h) = \sum_{k \ge 1}$ $\beta_t^k
\lan e_k , h \ran$ ($(e_k)$ denotes an orthonormal basis in $U$).
We will often use the nonnegative bounded operator $Q= B B^*$ on
$H$, where $B^*$ denotes the adjoint of $B$.

 We will make the following two basic assumptions:
\begin{Hypothesis} \label{ou1}
(i) \ the linear bounded operators $Q_t : H \to H$,
\beq \la{qt1}
Q_t x = \int_0 ^t e^{sA}
  Q e^{sA^*}x ds,\;  x \in H, \;\; \mbox { \rm   are trace class}, \;\; t>0;
\eeq
(ii) \  there exists  $ T>0, \;\; \mbox{\rm such that} \;\; e^{tA} (H)
\subset Q^{1/2}_t (H),\;\; t \ge T$.
\end{Hypothesis}
Assumtion (i) allows to solve equation (\ref{1e4}), by means of
the stochastic convolution, i.e. the Ornstein-Uhlebeck process
$X_t$ is given by
$$
X_t^x = e^{tA} x + \int_0^t e^{ (t-s) A } B d W_s,\;\; t \ge 0,\;
x \in H.
$$
Remark that if $H$, $U$ are finite dimensional, then Hypothesis
\ref{ou1} reduces to (\ref{qtt}), see for instance \cite{Za}.

Let us denote by $U_t$,  the transition Markov semigroup
associated to $X_t$,
 i.e. the Ornstein-Uhlenbeck semigroup.
 The following explicit formula, which generalises a
  formula due to Kolmogorov, holds:
\beq \label{calo}
U_t f(x)= \E f(X_t^x) \, =\, \int_{H}  f(e^{t A}x \,+ \, y )  \; {
N}(0 , Q_t ) dy,\;
 f \in {\cal B}_b (H  ),
\eeq
 $ x \in H,\; t>0$. Here $N(0, Q_t)$ denotes the  Gaussian measure on $H$ with mean 0
   and covariance operator $Q_t$.
 In \cite{DZ1} it is proved that
  the previous hypothesis (ii) is equivalent to the fact that
 $X_t$ is strong Feller and irreducible. More precisely,
  in \cite{DZ1} it is proved that, under (i)Hypothesis \ref{ou1},
\beq \label{eq}
 U_t ({\cal B}_b (H)) \, \subset \, {\cal C}_b ^{\infty} (H),\;\;
 t \ge T.
\eeq
 This fact implies that any BHF
$\psi$ for $U_t$ belongs in fact to ${\cal C}_b ^{\infty}(H)$.

  By  (ii),  the closed operator
 $Q_t^{-1/2} e^{t A}$ is a bounded operator on $H$, for any $t \ge T$.
  Moreover  assumptions (i) and (ii) implies that the semigroup
 $e^{tA}$ is compact for any $t \ge T$. To see this we write
 $ e^{TA}= Q_T^{1/2} Q_T^{-1/2} e^{T A} $ and note that by (ii)
  the operator $Q_T^{1/2}$ is
 Hilbert-Schmidt.

 The operator $Q_t^{-1/2}
 e^{t A}$  has also a clear control theoretic meaning. It is the minimal energy operator
 associated to the following the deterministic control system:
\beq \label{t1}
 \frac{dy}{dt} =  Ay(t) + Bu(t),\;\;\; y(0)= x \in H, \; t \ge 0,
 \eeq
where  $u : ( 0, + \infty) \to   U$ indicates the  control on the
system (\ref{t1}), see for instance \cite{Za}, \cite{CZ}. In the
proof of the Liouville theorem, we will need a result on  the
behaviour of the minimal energy operator as $t $ tends to $+
\infty$, see \cite{PZ}. Note that it is well known  that the map:
$t \mapsto$ $\|{Q_t}^{-1/2 } e^{t A }x \|$ is not increasing on
$[T, \infty)$, for any $x \in H$. The following statement  is a
particular case of a result proved in \cite{PZ}:

 \begin{Theorem} \la{PZ}   Assume Hypothesis \ref{ou1} holds.
  Then one has:
\beq \label{min}
 \lim_{t \to \ty }Q_t^{-1/2} e^{tA} x = 0, \;\; x \in
H, \;\; \mbox{\rm if and only if } \;\; s(A) = \sup \{ Re
(\lambda)\; :\; \lambda \in \sigma (A) \} \le 0.
\eeq
\end{Theorem}
Recall that since $e^{TA}$ is compact, then $\sigma(A)$ is
discrete and consists of eigenvalues of finite algebraic
multiplicity. We are going to prove the following result.
 \begin{Theorem} \label{m1}
 Assume Hypothesis \ref{ou1}. Then there exists a nonconstant BHF $\psi$ for
 $X_t$ if and only if $s(A) >0$.
\end{Theorem}
\underline{\it The necessity part}. \ \ Let us assume that $s(A)
\le 0$. Take any bounded harmonic function $\psi $. We need to
show that $\psi$ is constant.

First, applying the Cameron-Martin formula and hypothesis (ii), it
is not difficult to show that, for any $f \in {\cal B}_b (H)$,
 \beq \la{we}
<D U_t f (x), h>  =
 \int_{H }  f (e^{t A}x \, + y)  \lan Q_t ^{-1/2 } y,
Q_t^{-1/2 } e^{t A } h \ran \, N(0 ,Q_t)dy,
 \eeq
 $t \ge T$, $ h, \, x \in H$.  Here the map $y \mapsto $ $\lan Q_t ^{-1/2 } y,
k \ran  $ is a  Gaussian random variable with mean 0 and covariance $|k|^2$ defined
on the probability space $(H, N(0, Q_t))$, $k \in H$, $t \ge T$,  see \cite{DZ1},
Zabczyk \cite{KRZ} for more details.  By the H\"older inequality,  it follows that
\begin{eqnarray} \label{grad}
| \lan D U_t f (x),h \ran | \le \| f\| _0 \int_{H }  |  \lan  Q_t^{-1/2 } y,
Q_t^{-1/2 } e^{t A } h  \ran | \, N(0 ,Q_t)dy
\\
\nonumber \le |Q_t^{-1/2 } e^{t A } h | \| f\| _0, \;\; x \in H, \; t \ge T.
\end{eqnarray}
Then we argue similarly to the proof of Proposition \ref{grad}. We differentiate
both sides of the  identity $U_t \psi(x) = \psi(x)$. We get, using (\ref{grad}),
 \beq \la{we}
  \| \lan D  \psi (\cdot ), h \ran \|_0  \le
   |Q_t^{-1/2 } e^{t A } h | \| \psi \| _0,  \; t \ge T,\; h \in H.
\eeq
 Now letting $t \to \ty$,  we find that $\lan D  \psi (\cdot ), h \ran$
 is identically 0,  by Theorem \ref{PZ}. Since $h$ is arbitrary, we find that
 $\psi $ is constant and this proves the claim. \qed

In the proof of the sufficiency part,
  we need the following lemma which concerns  suitable prime integrals
 for linear deterministic dynamic systems.
\begin{Lemma} \label{dr2}
Let $A_0 \in {\cal L}(\R^n)$. Then there exists a
%%%%%% nonconstant, Lebesgue integrable,
Borel and bounded map $u : \R^n \to \R $, $\delta = \delta (A_0)>0$ such that, for
any $t \in [-\delta , \delta]$, one has:
  \beq \label{prime}
  u( e^{tA_0 } x) = u(x), \;\; x\not \in E_t,
  \eeq
 where $E_t$ has Lebesgue measure 0. Moreover, for any Borel probability measure
 $\mu$, one has that the map:
 \beq \label{vi}
 x \mapsto  \int_{\R^n} u( x + y ) \mu (dy),\;\;\; x \in \R^n,
\eeq
 differs from a constant  almost surely.
 \end{Lemma}
{\bf Proof} \ By means of a linear change of coordinates,
 we may assume that $A_0$ is in the canonical real Jordan form. Thus we consider
 coordinates $(x_1, \ldots, x_n )$ with respect to a Jordan basis of generalized
 eigenvectors of $A_0$.

First  let us consider the case in which there exists a real eigenvalue $\mu $ for
$A_0$. Let $x_k$ be the last variable in the Jordan block corresponding to $\mu$.
 For any function $f : \R^n \to \R$, depending only on $x_k$, i.e. $f(x) = \hat f(x_k)$,
 it is clear that
$$
f(e^{tA_0} x) = \hat f( e^{t\mu} x_k),\;\;\; t \in \R, \; x \in \R^n.
$$
Now if we take $u(x) = I_{(0, +\infty)}(x_k) = \hat u (x_k) $, we get
$$ u(e^{tA_0}
x) = \hat u (x_k) = u (x), \;\; t \in \R,
$$
 for any $x \in \R^n$, such that $x_k \not = 0$. This way (\ref{prime}) is satisfies.
To verify (\ref{vi}), note that
$$
\phi (x)= \int_{\R^n} u( x + y ) \mu (dy) = \mu (F_{x_k} ),\;\; x \in \R^n,
$$
where $F_{x_k} = \{ y \in \R^n \; : \; y_k >  - x_k  \}$, $x_k \in \R$. Arguing by
contradiction, if $\phi $ is equal to a constant $c$, a.s., then $c =1$. Now it is
easy to get a contradiction, by  using the complementary sets of $F_{x_k}$.

Let us consider the more difficult case  in which there exists a complex
 eigenvalue $\mu = a + ib $, $b \not =0$. It is not restrictive to assume that
  $A_0$ coincide with  the real Jordan block
 corresponding to $\mu$, i.e.
 $$
A_0=\left(
\begin{array}{cccccccc}
a &   b     & 1 & 0  & 0 & 0 & 0 &0    \\
-b &  a     & 0  & 1  & \cdots &\cdots & \cdots  & \vdots   \\
\vdots  & \vdots  & \ddots & \ddots & \ddots &  \ddots & \ddots & \vdots          \\
0       & \cdots  & 0  & \ddots & \ddots &  \ddots & 1 & 0 \\
0       & \cdots  & 0  & \ddots & \ddots & \ddots & 0 & 1  \\
0       & \cdots  & 0  & \ddots &  \ddots &\ddots  & a & b \\
0       & \cdots  & 0  & \ddots & \ddots & \ddots & -b & a  \\
\end{array}
\right)
$$
 If a  map $f$ depends only on  the last two variables of the
 previous Jordan block, i.e. $f (x)= \hat f (x_{n-1}, x_n)$, then
 $f (e^{tA_0} x) = \hat f ( e^{tB} (x_{n-1}, x_n))$,
$$
B = \left(
\begin{array}{cc}
a &       b        \\
-b       & a  \\
\end{array}
\right)
$$
By elementary methods,  we solve:
$$
\left \{ {\begin{array}{l}\displaystyle {
  x_{n-1}' (t) = a x_{n-1} + b x_n ,}
\\ \displaystyle{ x_{n}' (t) = a x_{n} - b x_{n -1} .}
\end{array}}\right.
$$
Eliminating the parameter $t$, we get a solution $x_n = x_n (x_{n-1})$, defined
implicitly by
\begin{eqnarray*}
\hat v (x_{n-1}, x_n (x_{n-1})) = \mbox {\rm constant},
\\
\mbox {\rm where} \; \hat v  ( x_{n-1}, x_n ) =
 \lg (\sqrt { x_{n-1}^2 +  x_n^2}) \, + \, \frac {a}{b} \arctan (\frac{
 x_n } {  x_{n-1} } ),  \; x_{n-1}  \not = 0,
 \end{eqnarray*}
$\hat v  ( 0, x_n ) = 0$. Now a direct computation shows that, for any $t \in
(-\delta , \delta)$, with $\delta = \frac {\pi }{ 2 |b|}$,
\begin{eqnarray*}
\hat v (e^{tB} (x_{n-1}, x_n) ) = \hat v (x_{n-1}, x_n), \;\;\; (x_{n-1}, x_n ) \not
\in E_t,
\end{eqnarray*}
where $E_t = \{ (x_{n-1}, x_n) \in \R^2 \; : \, x_{n-1}\cdot x_n \not = 0,\; x_{n-1}
\not = - \tan (tb) x_n   \}$.

Let us  introduce  the map: \beq \la{er3}
 u(x) = \exp[-(\hat v  ( x_{n-1}, x_n )   )^2  ] =
 \exp[-\Big( \lg (\sqrt { x_{n-1}^2 +  x_n^2})  +
  \frac {a}{b} \arctan (\frac{  x_n } {  x_{n-1} } )
 \Big)^2  ],
\eeq
 $x_{n-1} \not = 0$. We claim that $u$ verifies the assertions.
 To this end, let us check that it is Lebesgue integral on $\R^n$.  There results:
\begin{eqnarray*}
\int_{\R^n} u(x) dx =  \int_{\R^2}
 \exp[-( \lg (\sqrt { x_{n-1}^2 +  x_n^2})  +
  \frac {a}{b} \arctan (\frac{  x_n } {  x_{n-1} } )
 )^2  ] dx_{n-1}dx_{n}
\\
= \int_0^{2\pi} d\theta \int_{0}^{+\infty}
 \exp[-( \lg \rho  \, + \,
  \frac {a}{b} \arctan (\frac{  \sin \theta  } { \cos \theta } )
 )^2  ] \rho \, d\rho
\\
= \int_0^{2\pi} d\theta \int_{- \infty}^{+\infty}
 \exp[-( \lg (e^t)  \, + \,
  \frac {a}{b} \arctan (\frac{  \sin \theta  } { \cos \theta } )
 )^2  ] e^t \, dt
\\
\le {2\pi} \int_{- \infty}^{+\infty}
 \exp (- t^2)  \exp( t(1 + \frac {|a|}{|b|} {\pi \over 2} ) ) \; dt < + \infty.
\end{eqnarray*}
To prove (\ref{vi}), assume by contradiction that there exists $c \ge 0$, such that
$$
\int_{\R^n} u( x +  y ) \mu (dy) = c,\;\;  \; x \in \R^n,\; \mbox {\rm a.s.}.
$$
Integrating both sides with respect to the Lebesgue measure $\lambda $ and using the
Fubini Theorem, we get:
$$
\int_{\R^n} \mu (dy) \int _{\R^n} u( x +  y ) \lambda (dx) = \|u\|_{L^1} \; = \;
\int_{\R^n} c  \lambda (dx),
$$
where $L^1= L^1(\R^n, \lambda)$. This gives a contradiction since $\|u\|_{L^1}
>0$. The proof is complete. \qed

\noindent
 \underline{\it The sufficiency part}. Here we assume that $s(A)
 >0$ and construct a nonconstant bounded harmonic function $h$.
 Note that in order that $h$ is a BHF for
 $U_t$, it is enough that $U_t h = h$,  $t \in [0, \delta ]$, for some $\delta
 >0$. This follows by the semigroup property.

  We will  use a suitable decomposition of $H$.
 Since $e^{tA} $ is compact, $t \ge T,$ we can take an isolated eigenvalue
 $\mu  \in \sigma (A)$, Re$(\mu) >0$
 and we define $D_0$ as the finite dimensional
 subspace of all generalized eigenvectors associated to $\mu$.
  We introduce the  linear Riesz projection $P_0$
 (not orthogonal in general), see for instance pages 72 and 98 in \cite{CZ},
\beq \la{p0}
 P_0 x = \frac{1}{2 \pi i} \int_{\gamma} (w - A)^{-1}x dw,\;\;\; x \in H,
 \eeq
 where $\gamma$ is a circle enclosing $\mu$ in its interior
and $\sigma (A)/ \{ \mu \}$ in its exterior. We  get  $P_0 (H)= D_0$. Note that $H=
D_0 \oplus D_1$, where $D_1 = P_1 H$ and $P_1 = I- P_0$.  The closed subspaces $D_0$
and $D_1$ are both invariant for $e^{tA}$ and  moreover $ D_0 \subset D(A)$. We set
 $A_0= AP_0$, and
 \beq\label {ce3}
A_0 : D_0 \to D_0,\;\;\; A_1  : (D(A) \cap D_1) \subset D_1 \to
D_1.
 \eeq
 The operator $A_0$ generates a  group
 $e^{t A_0}$ on $D_0$ and $A_1 $
 generates a ${\cal C}_0$-semigroup $e^{tA_1}$ on $D_1$. The
 restrictions of $e^{tA}$ to $D_0$ and $D_1$ coincide with $ e^{t A_0}$ and
 $e^{tA_1}$ respectively. Moreover $A_0$ has the  unique eigenvalue $\mu $ on $D_0$.

 Let us assume that dim$(D_0)= n$ and   identify $D_0$ with $\R^n$.
 We consider the prime integral $u$ for the dynamical system $x' = A_0 x$,
  constructed in Lemma \ref{dr2}.

Let us introduce the operator:
$$
\ttq = \int_0 ^{\ty} e^{- s A_0} P_0 Q P_0 e^{- s A_0^* } x ds,\;\; x \in H,
$$
on $H$ (note that $\ttq (D_0) \subset D_0$) and define the Borel map
 \beq \label{inv} h(x) =
\int_{H} u(P_0 x + P_0 y ) N(0, {\ttq}) (dy), \;\;\; x \in H,
 \eeq where $u $ is
 given in (\ref{er3}). Note that $N(0, {\ttq}) $ is the invariant measure for the
 ``reverse'' Ornstein-Uhlenbeck process
$$
 \tilde X_t ^x = e^{- tA_0}x  + \int_0^t e^{-(t-s)A_0} B dW_s, \;\; x \in D_0, \; t \ge 0,
$$
 on $D_0$. By Lemma \ref{dr2}, we know that $h$
 is nonconstant, a.s. We are going to verify that  $h $ is  {\it harmonic} for $U_t$.
  Note that $P_0 e^{tA}P_0 = e^{tA_0}$ and $P_0^* e^{tA^*}P_0^* = e^{tA_0^*}$
 (the adjoint operations are all considered in $H$). We write, for any $t \in [0 ,
 \delta),$
\begin{eqnarray*}
 U_t  h (x)= \int_H  \int_{H} u ( P_0 e^{tA } x + P_0 z   + P_0 y  )  N _{  0 ,Q_t }  (dz)  N_{0 , \ttq }(dy)
\\
=\int_H  \int_{H} u ( e^{tA_0} [   P_0x +  e^{-tA_0}P_0 z   + e^{-tA_0}P_0 y ] )  N _{  0 ,Q_t }  (dz)  N_{0 , \ttq }(dy)
\\
=\int_H  \int_{H} u (  P_0 x +  e^{-tA_0}P_0 z   + e^{-tA_0}P_0 y  )  N _{  0 ,Q_t }  (dz)  N_{0 , \ttq }(dy)
\\
= \int_H  \int_{H} u (  P_0 x +  P_0 z+ P_0 y  )  N (  0 , e^{-tA_0} Q_t e^{-t[A_0]^*})  (dz)  N(0 , e^{-tA_0} \ttq e^{-t [A_0]^*})(dy)
\\
\int_H  u ( P_0 x +   P_0y  )  N  (  0 , e^{-tA_0} Q_t e^{-t[A_0]^*} +
 e^{-tA_0} \ttq e^{-t [A_0]^*} ) (dy),\;\; x \in H.
\end{eqnarray*}
Remark that
\begin{eqnarray*}
e^{-tA_0} Q_t e^{-t[A_0]^*}x +   e^{-tA_0} \ttq e^{-t [A_0]^*}x
\\
e^{-tA_0}P_0 (\int_0^t e^{sA} Q e^{sA^*}x ds ) P_0^* e^{-tA_0^*}
+ e^{-tA_0}( \int_0^{\ty} e^{-sA_0} P_0QP_0^* e^{- sA_0^*}x ds )  e^{-t [A_0]^*}
\\
= \ttq x,
\end{eqnarray*}
taking into account that $P_0$ is a projection and commutes with $e^{tA}$.
We get that $h$ is a bounded harmonic function for $U_t$.
\qed

Note that in order to prove the sufficiency part of Theorem \ref{PZ}, i.e. to
construct BHF for $U_t$, we need {\it less} than Hypothesis (ii) in (\ref{qt1}). We
only need that there exists an isolated $\mu \in \sigma (A)$, with Re$(\mu) >0$,
such that its associated invariant subspace of $A$ is finite dimensional.

\begin{Remark} \label{const}
{\em The previous Liouville theorem, with the same assumptions,
holds more generally the $H$-valued process $Y_t^x$, which solves:
\begin{eqnarray}   \label{gendeb}
 dY_t &= & AY_t dt  + Ak dt + B dW_t, \;\; \mbox {\rm i.e.}
\\
&& \nonumber Y_t ^x = e^{tA} x +  \int_0^t e^{s A} k ds + \int_0^t
e^{ (t-s) A } B d W_s,\;\; t \ge 0,\; x \in H,
\end{eqnarray}
where $k \in D(A)$ is fixed. To see this first note that the
process $X_t^{x+k} = Y_t^x + k $, verifies the SDE's $d X_t $ $ =
A (Y_t +  k) dt + B dW_t$ $= AX_t dt + B dW_t$. This means that
$X_t$ is an Ornstein-Uhlenbeck process, starting from $x+k$. Hence
we have:
 $$
T_t f(x) = \E f(Y_t^x) = \E f( X^{x +k}_t \,  - k ), t \ge 0,\; x
\in
 H.
 $$
This shows that there is a one to one  correspondence  between the
bounded harmonic functions for  $X_t$ and the ones for $Y_t$.
Indeed take for instance any BHF $g$ for $Y_t$, then define the
map $h$, $h (y ) = g (y - k)$, $y \in H$. Now  $h $ is a BHF for
$X_t$, because
$$
\E h (X_t^z) = \E h ( Y_t^{z-k} + k ) = \E g( Y_t^ {z-k}) = g(z-k) = h(z), \; z \in
H,\, t \ge 0.  \qed
$$
   }
\end{Remark}

We finish the section by a   result which shows that nonnegative harmonic functions
for the Ornstein-Uhlenbeck semigroup, non necessarily bounded, are always convex if
$s(A) \le 0$.
\begin{Proposition} \la{ci}
(\footnote {  The authors are indebted to  prof. Kwapien for this result }) Let
$\psi : H \to \R$ be a nonnegative function such that $U_t \psi = \psi $, $t \ge 0$.
Assume that  s$(A) \le 0 $, where s$(A)$ is the spectral bound of $A$, see
(\ref{min}). Then $\psi$ is convex.
\end{Proposition}
The proof is based on the following lemma.
\begin{Lemma} \la{kw} \
 For any nonnegative function $f : H \to \R$, there results:
 \beq \la{har}
 U_t f (x+a ) + U_t f(x+a ) \ge 2 C_t(a) \, U_t f(x), \;\; x,\, a \in H,
 \, t > 0,
 \eeq
  where $0 < C_t(a) < 1$, $t >0$. Note that both sides can be $+
  \infty$.
\end{Lemma}
{\bf Proof} \ By the Cameron-Martin formula one has:
\begin{eqnarray*}
  U_t  f (x+a)= \int_H  f (e^{tA}x  + y  )  \frac{ dN
_{e^{tA} a ,Q_t }  } { dN_{0 ,Q_t} }(y) N_{0 ,Q_t}(dy)
\\
= \int_H  f (e^{tA}x  + y  )   \exp [-\frac{1}{2} |Q_t^{-1/2}
e^{tA } a|^2 + <Q_t^{-1/2} e^{tA} a , Q_t^{-1/2}y> ] N_{0
,Q_t}(dy).
\end{eqnarray*}
It follows that
\begin{eqnarray*}
\frac {1}{2} ( U_t f (x+a) +  U_t f  (x-a ))
\\
 =  e^{ -\frac{1}{2} |Q_t^{-1/2} e^{tA } a|^2 }\,
 \int_H  f (e^{tA}x  + y  )
  \frac{1}{2}
  \Big ( e^{ <Q_t^{-1/2} e^{tA} a , Q_t^{-1/2}y> } + e^{ - <Q_t^{-1/2} e^{tA} a , Q_t^{-1/2}y> } \Big) N_{0 ,Q_t}(dy);
\\
\ge  \exp [-\frac{1}{2} |Q_t^{-1/2} e^{tA } a|^2 ] \, \int_H  f
(e^{tA}x  + y  )    N_{0 ,Q_t}(dy)
\\
 = \, C_t (a)\,   U_t f
(x),\;\;\; \mbox {\rm where } \; \, C_t (a) = \exp [-\frac{1}{2}
|Q_t^{-1/2} e^{tA } a|^2 ]. \qed
 \end{eqnarray*}
{\bf Proof of Proposition \ref{ci}}. \   By the previous lemma, we have:
\begin{eqnarray*}
\frac {1}{2} (\psi(x+a)  + \psi (x-a) ) = \frac {1}{2}
(U_t\psi(x+a) + U_t \psi (x-a) )
\\
\ge  \exp [-\frac{1}{2} |Q_t^{-1/2} e^{tA } a|^2 ] U_t \psi (x) =
\exp [-\frac{1}{2} |Q_t^{-1/2} e^{tA } a|^2 ]  \psi (x).
 \end{eqnarray*}
Passing  to the limit as $t \to \infty$, we infer, since $s(A) \le
0 $,
 \beq \la{cuy}
 \frac {1}{2} (\psi(x+a)  + \psi (x-a) ) \ge   \psi
(x),\;\; x,\, a \in H.
 \eeq
It is standard to verify that this condition is equivalent to the convexity of
$\psi$, see for instance ??. \qed

One should  notice that if in the previous Corollary the function
$\psi $ is in addition bounded, than formula  (\ref{cuy}) provides
 another proof of the necessity part of    Theorem \ref{m1}.
  Indeed in this case, we may assume that
  that $1 - \psi $  is a nonnegative BHF  (otherwise replace
$\psi$ by $\frac {\psi } {\| \psi\|_0 } $).  Replacing $\psi$ with
$1 - \psi $ in  (\ref{cuy}), we obtain:
\begin{eqnarray*}
\frac {1}{2} (1 - \psi(x+a) +1 - \psi (x-a )) = 1 - \frac{1}{2}
(\psi(x+a) + \psi (x-a )) \ge 1- \psi (x)
 \end{eqnarray*}
It  follows that $\psi(x+a) +  \psi (x-a ) \le  2 \psi (x)$ and
so, by (\ref{cuy}),
$$ \psi(x+a) +  \psi (x-a ) = 2 \psi
(x),\;\;\; x \in H.$$
 Now any continuous function $\psi$
 which satisfies the previous identity is affine, i.e. $\psi (x)= \psi (0) +
\lan h,x \ran$ for some $h \in H$. It follows that $\psi $ is a constant.

\section { Liouville theorems for Generalised Mehler   processes }

 in Section 4,   we
   a prove a Liouville theorem for     infinite
 dimensional Ornstein-Uhlenbeck processes perturbed   by a Levy noise.
 These processes are also called generalized Mehler processes and recently have received
a lot of attention, see for instance \cite{Za2}, \cite{Ch},
\cite{FR}, \cite{RW} and references therein. On this subject it is
also useful to recall  a few known result,
 which is proved   in \cite{CD}, see also \cite{Fe}, page 382:
 let  $Z_t$ be any Levy process
 with values in $\R^n$; then the only BHF for $Z_t$ are the constant ones.

In this section we  extend  the previous Liouville theorem to perturbations of
Ornstein-Uhlenbeck processes by means of jump processes. These are also called
generalized Mehler processes,  see  \cite{sato}, \cite{Za2}, \cite{Ch}, \cite{BRS},
\cite{FR},  \cite{RW} and  references therein. Let us consider the SDE's
 \beq \label{ju}
 dX_t = AX_t dt +    B dW_t  + C dZ_t ,\;\; X_0 = x
 \in H, \; t \ge 0,
 \eeq
where $A: D(A)\subset H \to H$ generates a ${\cal C}_0$-semigroup
$e^{tA} $ on $H$, $B$ and $C \in {\cal L}(U,H)$, where $U$ is
another Hilbert space. Moreover $W_t$ is an $U-$valued {\it
cylindrical} Wiener process adapted with respect to a stochastic
basis $(\Omega, {\cal F}, ({\cal F}_t), \P)$ and $Z_t $ is an
$U$-valued Levy process defined on the same stochastic basis and
adapted to ${\cal F}_t$. This means that  $Z_t$ is an $U-$valued
process with  with stationary independent increments, cadlag,
continuous in probability and such that $Z_0 =0$. Since a Gaussian
component is already present in (\ref{ju}), we consider $Z_t$
having only a pure jump part, i.e.
$$
Z_t = at + \xi_t, \;\; t \ge 0,
$$
where $a \in U$, and $\xi_t$ is a jump process having Levy
spectral measure $M$. This means  that  the law of $\xi_t $ has
characteristic function  given by
 \beq \la{ft1}
 \E e^{i \lan \xi_t , s \ran }  = \exp (t \phi (s) ), \; \mbox {\rm where } \; \phi(s)=
 \int_U \Big( 1 - e^{i \lan s,y \ran }  - \frac { i \lan s,y \ran} {1 +
|y|^2 } \Big )  M (dy), \;\; s \in U.
 \eeq
 Recall that $M $ is a $\sigma$-finite measure on $U$, satisfying $M (\{ 0 \}
 )=0$ and $
 \int_U (|y|^2 \wedge 1 ) M (dy )$ $ < \infty.$ Set $\varphi (s) $ $=
 \phi(s)  + i \lan a, s \ran$, $s \in U$.
  Let $Y_t$ be the Ornstein-Uhlenbeck process, which solves $dY_t
 = AY_t dt +    B dW_t  ,$ $ Y_0 = x
 \in H,$ $ t \ge 0. $ We assume that $Y_t$ is {\it independent }
 on $Z_t$.

 Following \cite{Ch} and \cite{DZ}, one can show that there exists a unique mild solution to (\ref{ju}),
which is given by
\begin{eqnarray} \label{dec}
X_t^x = e^{tA} x + \int_0^t e^{(t-s) A} B dW_s + \int_0^t e^{(t-s)
A} C d Z_s
\\
\nonumber = Y_t^x + \eta_t, \; \mbox {\rm where} \; \eta_t =
\int_0^t e^{(t-s) A} C d Z_s .
\end{eqnarray}
Note that the last stochastic integral can be defined as a limit
in probability of the elementary processes
 $$ \Phi_{n, t} =
\sum_{k=0}^{n-1} e^{(k t/n) A} C [Z_{(k+1) t/n } - Z_{k t/n }].
$$
Throughout this section, we assume Hypothesis \ref{ou1} and further
 that
 \beq \la{ou2}
  \;\;\; \int_U ( (\log|y| \wedge 0 ) \, M (dy ) < \infty.
 \eeq
 Note that if dim$(U)$, dim$(H) < \infty$ and $e^{tA}$ is
 stable, then (\ref{ou2}) is equivalent to the existence
  of an invariant measure for $X_t$.

Now we consider the characteristic function, or briefly c.f., of
the law of $X_t$, see also \cite{FR}. Since the c.f. of the law of
each $\Phi_{n,t}$ is $\exp(\sum_{k=0}^{n-1} \varphi (e^{(k t/n)
A^*} C^*h ) k/t )$, it follows that the law of $\int_0^t e^{ (t-s)
A } C d Z_s$ has c.f.
$$
\exp( \int_0^t  \varphi (e^{s A^*}  C^*h ) ds ) ,\;\; t \ge 0,\;
h\in U,
$$
where $C^* \in {\cal L}(H,U)$ is the adjoint of $C$ and $\varphi =
i \lan a, \cdot \ran$ $  + \phi $, and $\phi$ is given in
(\ref{ft1}). Now the law $\mu_t$ of $X_t^x $ is still infinitely
divisible on $H$, with c.f.
 \begin{eqnarray} \la{mil1}
 \hat \mu_t (h) = \exp (\psi_t (h)), \; \mbox {\rm where}\;
\\
\nonumber \psi_t (h) = i \lan a_t , h \ran - \frac{1}{2} \lan Q_t
h,h \ran + \int_H \Big( 1 - e^{i \lan h,y \ran  }  - \frac { i
\lan h,y \ran} {1 + |y|^2 }     \Big )  M_t (dy),
\\
\nonumber a_t = e^{tA} x   +  \int_0^t e^{sA}C a ds +  \int_0^t
\int_U e^{sA} C y \Big( \frac { 1} {1 + |e^{sA}Cy|^2 } -
 \frac { 1} {1 + |y|^2 }  \Big )  M (dy) ds,
\\
\nonumber
 Q_t h = \int_0^t e^{sA} B B^* e^{s A^*} h ds, \; h \in H,
\;\;  M_t (G) = \int_0^t M ( [ e^{sA} C]^{-1} \, G) ds, \;\;t \ge
0,
\end{eqnarray}
where $G$ is any Borel subset of $H$ and $[ e^{sA} C]^{-1}G $ is
the counterimagine  of $G$.  Note that
$$
\int_H (|y|^2 \wedge 1 ) M_t (dy ) = \int_0^t\int_U (|e^{sA} Cx|^2
\wedge 1 ) M (dx )  < \infty.
$$
After these preliminaries, we  prove the following result, which
generalizes a similar statement given in \cite{RW}. We denote by
$P_t$ the transition Markov semigroup associated to $X_t^x$.

\begin{Proposition}\la{RW}
 Let us assume that Hypothesis \ref{ou1} holds. Then $X_t^x $ is
 irreducible and moreover $P_t ({\cal B}_b (H)) \subset
 C_b^{\infty} (H)$, $t>0$.
\end{Proposition}
{\bf Proof} \ We write according to (\ref{dec}), $X_t^x = Y_t^x +
\eta_t$. Remark that, by Hypothesis \ref{ou1}, the
Ornstein-Uhlenbeck process $Y_t^x$ is irreducible. Using that
$Y_t^x$ and $\eta_t$ are independent,  it is easy to check that
$X_t^x $ is irreducible as well.

Let us prove that second assertion. Take $f \in {\cal B}_b (H) $
and denote by $\nu_t$ and $N(0,Q_t)$ the laws of $\eta_t$ and
$Y_t^0$. Since the  law $\gamma_t$ of $X^0_t$ is $\nu_t *
N(0,Q_t)$, it  follows that
\begin{eqnarray*}
P_t f(x) = \int_H f( e^{tA}x + y  ) \gamma_t (dy)
\\
= \int_H \nu_t (dz) \int_H f( e^{tA}x + y + z  ) N(0,Q_t) (dy)
\\
= \int_H \nu_t (dz) \int_H f(  y + z  ) N(e^{tA}x,Q_t) (dy), \; t
\ge 0, \; x\in H.
\end{eqnarray*}
Using the Cameron-Martin formula we can differentiate $P_t f$ in
each direction $h \in H$ and  get, for any $ x \in H,$ $t>0$,
 \beq \label{erw}
<D P_t f (x), h>  = \int_{H} \nu_t (dz)
 \int_{H }  f (e^{t A}x \, + y + z)  < {Q_t }^{-1/2 } y,
{Q_t}^{-1/2 } e^{t A } h >  N(0 ,Q_t)dy.
 \eeq
Similar formulas can be easily established for higher order
derivatives of $P_tf$. Now it is straightforward to verify that
$P_t f \in {\cal C}_b^{\infty} (H)$, $t>0$. This concludes the
proof.
 \qed
Now we prove a Liouville theorem for $P_t$.
\begin{Theorem} \label{m11}
 Let us assume that Hypothesis \ref{ou1} and (\ref{ou2}) hold.
 Then there exists a nonconstant BHF $\psi$ for $P_t$ if and only if $s(A) >0$.
\end{Theorem}
{\bf Proof } \ \underline{\it The necessity part}.  \ The proof is
similar to the one used in the necessity part of Theorem \ref{m1}.
 Let us assume that $s(A) \le 0$. Take any BHF
 $\psi $, i.e. $P_t \psi = \psi$, $t \ge 0$. We show that
$\psi$ is constant.  By (\ref{erw}), we get the estimate:
\begin{eqnarray*}
\| <D\psi(\cdot),h> \|_0 = \| <D P_t \psi( \cdot),h> \|_0
\\
\le \| \psi\| _0  \int_H \nu_t (dz) \int_{H }  |  < {Q_t }^{-1/2 }
y, {Q_t}^{-1/2 } e^{t A } h >| \, N(0 ,Q_t)dy
\\
\le  |{Q_t}^{-1/2 } e^{t A }h |  \| \psi\| _0,\;\; t>0,\, h \in H.
\end{eqnarray*}
Now letting $t \to \ty$ in the last formula,  we get that $\psi$
is constant, using that $|{Q_t}^{-1/2 } e^{t A } h | \to 0$ as $t
\to \infty $, see  Theorem \ref{PZ}. This proves the necessity
part.

\par
\hh \underline{\it The sufficiency part}  \ The proof
 of this part differs from   the corresponding one, considered in
 the sufficiency part of Theorem \ref{m1}.
 Here we will use  probabilistic arguments.

We assume that $s(A) >0$ and construct a nonconstant BHF $h$ for $P_t$. First we
take an isolated eigenvalue
 $\mu  \in \sigma (A)$, with Re$(\mu) >0$. Then
 we use the same suitable decomposition of $H$ used in the proof of Theorem \ref{PZ}
  with the same notation.  Let us assume that dim$(D_0)= n$.
 % When  no confusion may arise,
 %% we  will  identify $D_0$ with $\R^n$.
Suppose that we find $f :
D_0 \to \R$, such that
$$
\E f(P_0 X_t ^{a}) = f(a),\;\; a \in D_0.
$$
We claim that $h(x) = f(P_0 x )$, $x \in H$, is a BHF for $X_t^x$.
To see this we compute, taking into account  that $P_0$ commute
with $e^{tA}$,
 for any $x \in H,$
\begin{eqnarray*}
\E h(X_t^x) =  \E f (P_0 X_t ^{ x })
\\
= \E f (P_0 e^{tA} P_0 x + P_0 \int_0^t e^{(t-s) A} B dW_s +
\int_0^t e^{(t-s) A} P_0 C dZ_s)
\\
= \E f(P_0 X_t ^{P_0 x})=  f (P_0 x) = h(x).
\end{eqnarray*}
Remark that the previous computation shows that the  process
$\tilde X_t^a = $ $P_0 X_t^a $, $a \in D_0$, with values in $D_0$,
is given (up to equivalence in law) by
$$
\tilde X_t^a = e^{tA}a +  \int_0^t e^{(t-s) A} P_0 B P_0  d \tilde
W_s  + \int_0^t e^{(t-s) A} P_0 C dZ_s
$$
(this can be checked by using  characteristic function) where $\tilde W_t$  is a
standard $n$-dimensional Wiener process with values in $D_0$. Let us define the Levy
process $\tilde L_t = P_0 C Z_t $  $+ P_0 B P_0  d \tilde W_t $ with values in
$D_0$. We can then write:
 \beq \label{ci1}
\tilde X_t^a = e^{tA}a +  \int_0^t e^{(t-s) A} dL_t.
 \eeq
 Note that   the spectral measure $\tilde M$ of $L_t$,
 which  is given by $\tilde M (G)$ $= M ( (P_0 C)^{-1} G )$, for any
 Borel set $G$ in $\R^n$, satisfies
\begin{eqnarray*}
\int_{\R^n}  (\log|y| \wedge 0 ) \, \tilde M (dy ) = \int_{U} (
\log| P_0 C y| \wedge 0 ) \, M (dy )  < \infty,
\end{eqnarray*}
by assumptions (\ref{ou2}).
 Now we finish the proof by using the following
 lemma.

 \begin{Lemma} \la{costr}
 Let $D \in {\cal L } (\R^n)$ such that all its eigenvalues have  positive real part. Let
$$
 Y_t^x = e^{tD}x +  \int_0^t e^{(t-s) D}  d L_s , \;\; t \ge 0, \, x \in
 \R^n
$$
 where $L_t$ is any  Levy process in $\R^n$. Assume that
 \beq \la{ou3}
  \int_{\R^n} ( \log|y| \wedge 0 ) \, M (dy ) < \infty,
 \eeq
 where $M$ is the spectral measure of $L_t$. Then there exists a nonconstant BHF $f$ for
 $Y_t$.
\end{Lemma}
{\bf Proof} \ Let us consider the ``reverse process'' $ \hat X_t$, corresponding to
the matrix $-D$, i.e.
$$
\hat X_t^x =  e^{ - t D}x +  \int_0^t e^{ (s -t) D}  d  L_s .
$$
This process admits an invariant measure $\mu$, since $-D$ is
stable. Take a random variable $Z$ such that its law is $\mu$ and
further is independent on $X_t^x$, $t \ge 0$. One has:
 \beq \la{conv}
 \E g (  e^{-tA}Z + \int_0^t e^{ (s-t) D}  d  L_s  ) = \E g(Z), \;\;
 t \ge 0, \; g \in {\cal C}_b (\R^n).
 \eeq
 As in the proof of Theorem \ref{m1}, we consider the prime integral
 $u : \R^n \to \R$, constructed in Lemma \ref{dr2}.  Recall that $u$ verifies:
 \beq \la{dr}
 u( e^{t D } x) = u(x), \; x \in \R^n ,  \; t \in (-\delta , \delta),
 \eeq
almost surely with  respect to the Lebesgue measure. Now let us
define
$$
f(x) = \E u(x + Z),\;\; x \in \R^n.
$$
By Lemma \ref{dr2}, we know  that $f$ is not constant, a.s.. Now we check that $f$
is a BHF for $Y_t^x$. By the independence of $Z$ from $Y_t^x$, we obtain:
\begin{eqnarray*}
P_t f(x) = \E u (Y_t^x  + Z) = \E u(  e^{tD}x +  \int_0^t e^{(t-s)
D}  d  L_s  \, + \, Z )
 \\
= \E u(  x +  \int_0^t e^{-s D}   dL_s  \, + \, e^{-tA}Z), \;\; t \in [0, \delta)
\end{eqnarray*}
by (\ref{dr}). Now note that the laws  of $ \int_0^t e^{-s D} dL_s
$ and $\int_0^t e^{-(t-s) D}  d L_s  $  are the same, for any
$t\ge 0$. Indeed   their characteristic functions coincide with
$$
 \exp \Big( - \int_0^t  \psi (e^{sD^*} r) ds \Big ), \;\; r \in
 \R^n, \;\;  [0, \delta),
$$
where $\E e^{i \lan L_1 , r \ran }$ $=  e^ {\psi(r)} $. Using this
fact and (\ref{conv}) we infer
$$
P_t f(x) = \E u ( x+  e^{-tA}Z +  \int_0^t e^{ (s-t) D}  d  L_s )
= \E u(x+ Z) = f(x),
$$
for any $t  [0, \delta)$. This shows that $f$ is a nonconstant BHF for $Y_t$ and
concludes the proof. \qed

Remark that Proposition \ref{ci} and Lemma \ref{kw} hold even for the Generalized
Mehler process $X_t$, which solves (\ref{ju}), under Hypothesis \ref{ou1} and
(\ref{ou2}).

\begin{Remark} \label{FUR}
{\em  Let us mention that the approach in \cite {FR} is slight
different from our. Indeed they start from the Mehler semigroup
$$
P_t f(x) = \int_{H} f(e^{tA} x + y ) \mu_t (dy),
$$
and then they construct, under suitable assumtions, a process on a bigger space
having $P_t$ has corresponding transition Markov semigroup. We can also formulate
our Theorem \ref{m11} also in this setting. { \bf what do you think ?? Perhaps one
can give more details? }  I could also write a construction of BHF by means of
Fourier transform, in the spirit of Fuhrmann and Rockner. Or  could I put this proof
in Appendix ??   }
\end{Remark}

\section { Probabilistic interpretation of the Liouville theorem }

Here  we investigate in some generality the connection  between
existence of BHF $u$ for the transition Markov semigroup $P_t$,
see (\ref{he}), and the path properties of the corresponding
Markov process $Y_t$. We will assume that the process $Y_t$,  is
strong Feller and irreducible, and admits c\`adl\`ag trajectories.

%%%%%%%%%%%%5
Our main references are Chapter 9 of the the book \cite{Pi}, which
investigates in a deep way
  nondegenerate diffusions,  associated to strongly elliptic operators, in finite dimension.
 In Chapter 9 of \cite{Pi}, BHF for some nondegenerate diffusions are studied.
 This is done by combining analytic and probabilistic techniques.
  On this subject,  we also mention  \cite{CD},  \cite{Cr} and the
   references in \cite{Pi}.
%%%%%%%%%%%%%%%

\subsection{Formulation of the main result}

Let $(E, \rho)$ be a Polish space, i.e. a complete separable
metric space, with metric $\rho $ and the $\sigma-$algebra of all
Borel subsets indicated by ${\cal E}$. Let $P_t$ be a Markovian
semigroup determined by a transition function $p$, see
(\ref{mark}), and acting on ${\cal B}_b(E)$.

As before a map $u \in {\cal B}_b (E)$ is called BHF for $P_t$ if
 \beq \label{he1}
   P_t u(x) =u(x),\;\;\; t \ge 0, \; x \in E,
   \eeq
 Note that in order that $\phi$ is a BHF for
 $P_t$, it is enough that (\ref{he1}) holds for $t \in [0, \delta ]$, for some $\delta
 >0$. Now we will make the following basic assumptions on $P_t$:
\begin{Hypothesis}\label{fi} {\em
 \hh
 (i)\ \ for any $x \in E $, $r >0$, there exists $t = t(x,r) >0$ such that
 $p (t,x, \{ y \, : \, \rho (x,y)  < r \} ) >0$;
 \hh
(ii) \ \ there exists $t>0$ such that for any $f \in {\cal B}_b (E)$,
 $P_t f $ is continuous on $E$;
\hh
 (iii) \  \  Markovian measures determined by $P_t$ are supported
in $D_E[0, \infty)$,  the Skorokhod space of all c\`adl\`ag
functions $\omega: [0, \infty) \to E$ (i.e. $\omega $ is  right
continuous on $[0, \infty)$ and has  left limits at each  point
$t>0$).
 }
\end{Hypothesis}
Let us comment the previous assumptions. Properties (i) and (ii)
are weak versions of the so called irreducibility and strong
Feller property for Markov processes, see \cite{DZ1}. Note that,
in several cases, assumption (ii) holds when the transition
probability functions  verify:
 }}

%\vspace{ 5 mm } {\small
% Dipartimento di Matematica ``Giuseppe Peano'', Universit\`a
%di Torino \par
%  via Carlo Alberto 10  \ 10123
%  \par  Torino, Italy \par
% e-mail: enrico.priola@unito.it }
%\par \ \par

\end{document}